%


\countdef\stylenumber=1 
\countdef\footnotenumber=2 
\newif\ifamstexsupported
     \amstexsupportedtrue
       
\def\smallprint{\magnification=\magstep0
     \overfullrule0pt
     \def\pagewidth{468pt}
     \def\pageheight{9truein}
     \hsize\pagewidth
     \vsize\pageheight
     \stylenumber=2}

\def\discrversionA#1#2{\ifnum\stylenumber=1#1\else#2\fi}
\def\discrversionB#1#2{\ifnum\stylenumber=2#1\else#2\fi}
\def\discrversionC#1#2{\ifnum\stylenumber=3#1\else#2\fi}
\def\discrversionD#1#2{\ifnum\stylenumber=4#1\else#2\fi}
\def\discrversionE#1#2{\ifnum\stylenumber=5#1\else#2\fi}
\def\discrversionF#1#2{\ifnum\stylenumber=6#1\else#2\fi}
\def\discrversionG#1#2{\ifnum\stylenumber=7#1\else#2\fi}
\def\discrversionH#1#2{\ifnum\stylenumber=8#1\else#2\fi}
\def\filename#1{\ifnum\stylenumber=1{\noindent filename #1}
               \else\ifnum\stylenumber=6{\noindent{\fiverm filename #1}}
               \else\ifnum\stylenumber=7{\noindent{\fiverm filename #1}}
               \else\fi\fi\fi}
\def\versiondate#1{
    \immediate\write0{Version of #1}
               \ifnum\stylenumber=1\hfill Version of #1 
                              \medskip
               \else\ifnum\stylenumber=2\hfill Version of #1
                              \medskip
               \else\ifnum\stylenumber=6\hfill{\fiverm Version of #1}
                              \medskip
               \else\ifnum\stylenumber=7\hfill{\fiverm Version of #1}
                              \medskip
               \else\ifnum\stylenumber=8\hfill Version of #1
                              \medskip
               \else\fi\fi\fi\fi\fi}
\def\discrpage{\ifnum\stylenumber=1{\frnewpage}
           \else\ifnum\stylenumber=3{\frnewpage}
           \else\ifnum\stylenumber=5{\frnewpage}
           \else\bigskip\fi\fi\fi}
\def\frnewpage{\vskip \pageheight plus 0pt minus\pageheight
           \eject
           \gdef\topparagraph{}
           \gdef\bottomparagraph{}
           }
\let\oldfootnote=\footnote
\def\query{\discrversionA{\immediate\write0{query}
    \advance\footnotenumber by 1
    \oldfootnote{$^{\the\footnotenumber}$}{query}}{}}
\def\versionref{\ifnum\stylenumber=1\workingref
           \else\ifnum\stylenumber=2\workingref
           \else\ifnum\stylenumber=8\workingref
           \else\printerref\fi\fi\fi}
\def\versionhead{\ifnum\stylenumber=3 \hfill PRINTER'S VERSION

                        \bigskip
                   \else\ifnum\stylenumber=4 \hfill REFEREE'S VERSION

                        \bigskip
                   \else \fi\fi}

\def\essexaddress#1{\hskip 3truein Mathematics Department\par
   \hskip 3truein University of Essex\par
   \hskip 3truein Colchester CO4 3SQ\par
   \hskip 3truein England\par
   \noindent e-mail: {\tt fremdh\@essex.ac.uk}\hskip 1.5truein#1
   \bigskip}

\def\shortaddress{\noindent Mathematics Department, University of Essex,   
   Colchester CO4 3SQ, England\par
   \noindent e-mail: fremdh\@essex.ac.uk}

\def\Bbbone{\hbox{1\hskip0.13em\llap{1}}}

\def\frsmallcirc{{\hskip 1pt plus 0pt minus 0pt}
     {\raise 1pt\hbox{$\smallcirc$}}{\hskip 1pt plus 0pt minus 0pt}}

\long\def\leaveitout#1{}

\def\Loadfont#1{\font#1}

\def\restr{\hbox{$\hskip-0.05em\restriction\hskip-0.05em$}}

\def\smallcirc{{\scriptstyle\circ}}

\def\tbf#1{\text{\bf #1}}


\ifamstexsupported\else
   \immediate\write0{AmSTeX-unsupported version of  fremtex.tex }
   \def\@{@}
   \def\\{\cr}
   \def\Bbb#1{\underline{\hbox{\bf #1}}}

   \def\Cal#1{{\cal#1}}

   \def\frac#1#2{{#1}\over{#2}}
   \def\frak#1{{\hbox{\bf #1\hskip0.1em\llap #1}}}
   \def\frnewpage{\vfill\eject}

   \def\Loadfont#1{}

   \def\preccurlyeq{\preceq}
   \def\query{\discrversionA{\immediate\write0{query}
                             \footnote{$^*$}{query}}{}}
   \def\restr{|}

   \def\smallfrown{\wedge}
   \def\smc{}

   \def\sqr#1#2{{\vcenter{\vbox{\hrule height.#2pt 
          \hbox{\vrule width.#2pt height#1pt \kern#1pt
                  \vrule width.#2pt}
          \hrule height.#2pt}}}}
   
   \def\text#1{\hbox{#1}}

   \def\Vdash{|\vdash}
   
   \fi


\ifamstexsupported\input amstex
               \documentstyle{amsppt}

               \font\varseveneufm=eufm10 at 7 pt
               \scriptfont\eufmfam=\varseveneufm
   \def\=#1{{\accent"16 #1}}
   
   \else\fi

\filename{FSp90a.tex}
\versiondate{16.9.92}
\def\formset#1{\{#1\}}

\def\BJSpHI{[BJSp89]}
\def\BaHD{[Ba84]}
\def\FrpHI{[Frp89]}
\def\KuHJ{[Ku80]}
\def\KVHD{[KV84]}
\def\ShHB{[Sh82]}
\def\ShCBF{[Sh326]}
\def\SpHG{[Sp87]}
\def\TalHD{[Ta84]}
\def\TaHG{[Ta87]}

\smallprint

\centerline{\bf Pointwise compact and stable sets}

\centerline{\bf of measurable functions}

\vskip 12pt

\centerline{\smc S.Shelah \& D.H.Fremlin}

\centerline{\it Hebrew University, Jerusalem}

\centerline{\it University of Essex, Colchester, England}

\vskip 12pt

\discrversionA{
\centerline{[University of Essex Mathematics Department Research Report
91-3]}}{}
\discrversionB{
\centerline{[University of Essex Mathematics Department Research Report
91-3]}}{}
\discrversionE{
\centerline{[University of Essex Mathematics Department Research Report
91-3]}}{}

\discrversionE{\vfill\eject}{}

\vskip 12pt

{\bf Introduction} In a series of papers culminating in \TalHD,
M.Talagrand, the second author and others investigated at length the properties
and structure of pointwise compact sets of measurable functions.   A number of
problems, interesting in themselves and important for the theory of Pettis
integration, were solved subject to various special axioms.   
It was left unclear
just how far the special axioms were necessary.   In particular, several
results depended on the fact that it is consistent to suppose that
every countable relatively pointwise compact set of Lebesgue measurable
functions is `stable' in Talagrand's sense;  the point being that stable
sets are known to have a variety of properties not shared by all pointwise
compact sets.   In the present paper we present a model of set theory in
which there is a countable relatively pointwise compact set of Lebesgue
measurable functions which is not stable, and discuss the significance of
this model in relation to the original questions.   
A feature of our model which may be of independent interest is the following:
in it, there is a closed negligible set $Q\subseteq[0,1]^2$ such that whenever
$D\subseteq[0,1]$ has {\it outer} measure 1 then 

\centerline{$Q^{-1}[D]=\{x:\exists\enskip y\in D,\,(x,y)\in Q\}$}

\noindent has {\it inner} measure 1 (see 2G below).

\medskip

\noindent {\bf 1. The model}  We embark immediately on the central ideas of
this paper, setting out a construction of a partially ordered set which
forces a fairly technical proposition in measure theory (1S below);  the
relevance of this proposition to pointwise compact sets will be discussed
in $\S 2$.   The construction is complex, and rather than give it in a
single stretch we develop it cumulatively in 1E, 1I, 1Q 
below;  it is to be understood that each notation introduced in these
paragraphs, as well as those in the definitions 1A, 1K, 1L, 
is to stand for the remainder of the section.
After each part of the construction we give lemmas which can be
dealt with in terms of the construction so far, even if their motivation is
unlikely to be immediately clear.

When we come to results involving Forcing, we will try to follow the
methods  of \KuHJ;  in particular, in a p.o.set,
`$p \le q$' will always mean that $p$ is a stronger condition than $q$.

\medskip

{\bf 1A Definition} If $\Cal A$ is any family of sets not containing 
$\emptyset$, set

\centerline{dp$(\Cal A) = \min\formset
{ \#(I):I\cap A \ne \emptyset\enskip\forall
\enskip A \in \Cal A}$.}

\noindent Observe that dp$(\Cal A) = 0$ iff $\Cal A = \emptyset$ and that
dp$(\Cal A\cup\Cal B)$ is at most the cardinal sum of dp$(\Cal A)$ and
dp$(\Cal B)$.   (Of course much more can be said.)

\medskip

{\bf 1B Lemma} Suppose that $n$, $l$, $k\in\Bbb N$, with $n$, $l$
not less than $2$,
and that $\epsilon$ is such that $0<\epsilon\le\frac12$ and
$l\epsilon^k\ge (k+2)\ln n$.
Then there is a set $W\subseteq n\times n$ (we identify $n$ with the set of
its predecessors) such that $\#(W)\le\epsilon n^2$ 
and whenever $I\in[n]^l$ and $J_0,\ldots,J_{l-1}\in [n]^{\le k}$ are 
disjoint, 
there are $i\in I$, $j<l$ such that $\formset{i}\times J_j\subseteq W$.

\medskip

\noindent {\bf proof} If $k=0$ this is trivial;  suppose that $k>0$.
Set $\Omega=\Cal P(n\times n)$.
Give $\Omega$ a probability for which the events $(i,j)\in W$, as
$(i,j)$ runs over $n\times n$, are
independent with probability $\epsilon$.   If $W\in\Omega$ is a
random set, then

\centerline{Pr$(\#(W)\le\epsilon n^2) 
>\frac14$}

\noindent because $\epsilon\le\frac12$ and $\#(W)$ has the
binomial distribution
$B(n^2,\epsilon)$.
On the other hand, if $J\in [n]^{\le k}$ and $i<n$, 
Pr$(\formset{i}\times J
\subseteq W)\ge\epsilon^k$.   So if $I\in [n]^l$ and $J_0,\ldots,J_{l-1}$
are disjoint members of $[n]^{\le k}$,

$$\eqalign{\text{Pr}(\formset{i}\times J_j\not\subseteq W\enskip\forall\enskip
i\in I,\,j<l)&\le(1-\epsilon^k)^{l^2}\cr
&\le \exp(-l^2\epsilon^k).\cr}$$

\noindent Accordingly
\discrversionC
{ the probability that there are $I\in[n]^l$,
disjoint $J_0,\ldots,J_{l-1}\in
[n]^{\le k}$ such that 
$\formset{i}\times J_j\not\subseteq W\enskip\forall\enskip
i\in I,\,j<l$ is at most
}{}

\discrversionC
{$$\eqalign{\#([n]^l)\#([n]^{\le k})^l\exp(-l^2\epsilon^k)
&\le n^l n^{kl}\exp(-l^2\epsilon^k)\cr
&=\exp((k+1)l\ln n-l^2\epsilon^k)\le\frac14\cr}$$
}
{$$\eqalign{\text{Pr}(\exists\enskip I\in[n]^l,\,
\text{disjoint }&J_0,\ldots,J_{l-1}\in
[n]^{\le k}\text{ such that }\formset{i}\times J_j\not\subseteq W\enskip\forall\enskip
i\in I,\,j<l)\cr
&\le\#([n]^l)\#([n]^{\le k})^l\exp(-l^2\epsilon^k)\cr
&\le n^l n^{kl}\exp(-l^2\epsilon^k)\cr
&=\exp((k+1)l\ln n-l^2\epsilon^k)\le\frac14\cr}$$
}

\noindent because

\centerline{$l^2\epsilon^k-(k+1)l\ln n\ge l\ln n\ge 2\ln 2$.}

\noindent There must therefore be some $W\in\Omega$ of the type required.

\medskip

\noindent {\bf Remark} Compare the discussion of cliques in random graphs
in \SpHG, pp. 18-20.

\medskip

{\bf 1C Lemma} Let $m$ and $l$ be strictly positive integers and $\Cal A$
a non-empty family of non-empty sets.   Let $\Bbb T$ be the family of non-empty
sets $\Cal T\subseteq\Cal A^m$.   For $\Cal T\in\Bbb T$ write $\Cal T^*=
\formset{\bold t\restr j:\bold t\in\Cal T,\,j\le m}\subseteq\bigcup_{j\le m}
\Cal A^j$.   For $\Cal T$, $\Cal T_0\in\Bbb T$ say that $\Cal T
\preccurlyeq\Cal T_0$ if $\Cal T\subseteq\Cal T_0$ and

\centerline{dp$(\formset{u:\bold t^{\smallfrown}u\in\Cal T^*})\ge
\text{dp}(\formset{u:\bold t^{\smallfrown}u\in\Cal T_0^*})/2l$}

\noindent for every $\bold t\in\Cal T^*\setminus\Cal T$.   Fix
$\Cal T_0\in\Bbb T$ and a cover $\langle \Cal S_i\rangle_{i<2l}$ of
$\Cal T_0$.   Then there is a $\Cal T\preccurlyeq\Cal T_0$ such that
$\Cal T\subseteq\Cal S_i$ for some $i<2l$.

\medskip

\noindent {\bf [Notation:}  In this context we use ordinary italics, `$u$',
for members of $\Cal A$, and bold letters, `{\bf t}', for finite sequences
of members of $\Cal A$.]

\medskip

\noindent{\bf proof} For $\bold t\in\Cal T_0^*\setminus\Cal T_0$ set

\centerline{$\alpha_{\bold t}=\text{dp}(\formset{u:\bold t^{\smallfrown}u
\in\Cal T_0^*})/2l>0$.}

\noindent For $i<2l$ define $\langle\Cal S^{(j)}_i\rangle_{j\le m}$ by setting
$\Cal S_i^{(m)}=\Cal S_i\cap\Cal T_0$,

\centerline{$\Cal S_i^{(j)}=\formset{\bold t:\bold t\in\Cal A^j\cap\Cal T_0^*,\,
\text{dp}(\formset{u:\bold t^{\smallfrown}u\in\Cal S^{(j+1)}_i})
\ge\alpha_{\bold t}}$}

\noindent for $j<m$.   An easy downwards induction (using the fact that dp is
subadditive) shows that $\Cal T_0^*\cap\Cal A^j=\bigcup_{i<2l}\Cal S_i^{(j)}$
for every $j\le m$.   In particular, there is some $i<2l$ such that
$\emptyset\in\Cal S_i^{(0)}$.   Now define $\Cal T$ by

\centerline{$\Cal T=\formset{\bold t:\bold t\in\Cal A^m,\,\bold t\restr j
\in\Cal S^{(j)}_i\enskip\forall\enskip j\le m}\subseteq\Cal T_0\cap \Cal S_i$,}

\noindent and see that $\Cal T\preccurlyeq\Cal T_0$, as required.

\medskip

{\bf 1D Corollary} Let $n$, $l$, $k$ and $W$ be as in Lemma 1B.   Take $r\le k$,
let $Z$ be the cartesian product $n^r$ and set

\centerline{$\tilde W=\formset{(i,z):i<n,\,z\in Z,\,(i,z(j))\in W\enskip
\forall\enskip j<r}$.}

\noindent Let $m$, $\Cal A$, $\Bbb T$ and $\preccurlyeq$ be as in Lemma 1C,
and take $\Cal T_0\in\Bbb T$, $H:\Cal T_0\to n$ any function.   Then

{\it either} there are $i<n$, $\Cal T\preccurlyeq \Cal T_0$ such that
$H(\bold t)=i$ for every $\bold t\in\Cal T$

{\it or} there is a $J\in[n]^{\le rl}$ 
such that for every $z\in(n\setminus J)^r$
there is a $\Cal T\preccurlyeq\Cal T_0$ such that $(H(\bold t),z)\in\tilde W$
for every $\bold t\in\Cal T$.

\medskip

\noindent{\bf proof}
Set

\centerline{$A=\formset{z:z\in Z,\,\exists\enskip\Cal T\preccurlyeq\Cal T_0
\text{ such that }(H(\bold t),z)\in\tilde W\enskip\forall\enskip\bold t
\in\Cal T}$.}

\noindent If $A\supseteq(n\setminus J)^r$ for some $J\in[n]^{\le rl}$, we have
the second alternative;  suppose otherwise.   Then we can find $z_0,\ldots,
z_{l-1}\in Z\setminus A$ such that the sets $J_j=\formset{z_j(i):i<r}$ are all
disjoint.   Each $J_j$ belongs to $[n]^{\le k}$, so by the choice of
$W$,

\centerline{$I=\formset{i:\formset{i}\times J_j\not\subseteq W\enskip
\forall\enskip j<l}$}

\noindent has cardinal less than $l$.   Now observe that if $\bold t\in\Cal
T_0$ then either $H(\bold t)\in I$ or $(H(\bold t),z_j)\in\tilde W$
for some $j<l$.   So we have a cover of $\Cal T_0$ by the sets

\centerline{$\Cal S_i=\formset{\bold t:H(\bold t)=i}$ for $i\in I$,}

\centerline{$\Cal S_j'=\formset{\bold t:(H(\bold t),z_j)\in\tilde W}$
for $j<l$.}

\noindent By Lemma 1C, there is a $\Cal T\preccurlyeq\Cal T_0$ such that
either $\Cal T\subseteq\Cal S_i$ for some $i\in I$ or $\Cal T\subseteq
\Cal S_j'$ for some $j<l$.  But we cannot have $\Cal T\subseteq \Cal S_j'$,
because $z_j\notin A$;  so $\Cal T\subseteq\Cal S_i$ for some $i$, and
we have the first alternative.

\medskip

\noindent{\bf Remark} 1C-1D are of course elementary, but their significance
is bound to be obscure;  they will be used in 1R below.
An essential feature of 1C is the fact that the denominator `$2l$'
in the definition of $\preccurlyeq$ is independent of the size of $\Cal A$.

\medskip

{\bf 1E Construction:  part 1 (a)} Take a sequence $\langle n_k\rangle
_{k\in\Bbb N}$ of integers increasing so fast that

\quad(i) $n_0\ge 4$;

\quad(ii) $n_k> 2^{k+1}$;

\quad(iii) writing $\tilde c_l=\prod_{i<l}2^{n_i}$, then

\qquad$\ln(2^{-k-1}n_k)\ge 2^{l}(k+1)(\tilde c_l^{l+1}
\ln 2+\tilde c_l^{l-1}\ln n_{l-1})$ for $1\le l\le k$;

\quad(iv) $\ln n_k\ge (2^{k+1})^k(k+2)$;

\quad(v) writing $\lceil\alpha\rceil$ for the least integer greater than
or equal to $\alpha$,

\qquad$(k+1)\ln(2\lceil(\ln n_k)^2\rceil)
\le2^{-k}\lceil\ln(2^{-k-1}n_k)\rceil$;

\quad(vi) $2^kk\lceil (\ln n_k)^2\rceil(\prod_{i<k}2^{n_i})^{k+1}\le n_k $;   

\quad(vii) $\ln(2^{-k-1}n_k)\ge(k+1)\ln(2\tilde c_k^{k+1}+2k)$

\noindent for every $k\in\Bbb N$.   

For each $k\in\Bbb N$, let $V_k$ be the cartesian product $\prod_{i<k}n_i$.

\medskip

{\bf (b)} For each $k\in\Bbb N$, let $T_k$ be the set of those subsets
$t$ of $V_k$
expressible as $t=\prod_{i<k}C_i(t)$ where $C_i(t)\subseteq n_i$ and
$\#(C_i(t))\ge(1-2^{-i-1})n_i$ for each $i<k$.   Set $T=\bigcup_{k\in\Bbb N}
T_k$, and for $t\in T$ say that rank$(t)=k$ if $t\in T_k$.   For $t$, $t'
\in T$ say that $t\le t'$ if rank$(t)\le\text{rank}(t')$ and $C_i(t)=C_i(t')$
for every $i<\text{rank}(t)$.   Then $T$ is a finitely-branching tree of height
$\omega$ in which the $T_k$ are the levels and `rank' is the rank function.
For $t\in T$ write $T^{(t)}$ for the subtree $\formset{t':t'\le t\text{ or }
t\le t'}$, suc$(t)$ for $\formset{t':t\le t',\,\text{rank}(t')=\text
{rank}(t)+1}$.

\medskip

{\bf (c)} For $k\in \Bbb N$, set

\centerline{$\gamma_k=(k+1)/\ln(\lceil 2^{-k-1}n_k\rceil)$;}

\noindent 
$2^{-k-1}n_k>1$ by (a)(ii) above. 
For $t\in T$ define ${d}_t:\Cal P T\to\Bbb R\cup\formset{-\infty}$
by writing

\centerline{${d}_t(S) = \gamma_{\text{rank}(t)}\ln(\text{dp}
(\formset{C:t\times C\in S}))$}

\noindent for every $S\subseteq T$, allowing ${d}_t(S)=-\infty$
if $S\cap\text{suc}(t)=\emptyset$.   Observe that ${d}_t(T)\ge
k+1$ whenever $\text{rank}(t) = k$ (because

\centerline{dp$(\formset{C:C\subseteq n_k,\,\#(C)\ge(1-2^{-k-1})n_k})
\ge \lceil 2^{-k-1}n_k\rceil$.)}

\medskip

{\bf (d)} Let $\Bbb Q$ be the set of subtrees $q\subseteq T$ such that

\enskip $q\ne\emptyset$;

\enskip if $t\le t'\in q$ then $t\in q$;

\enskip if $t\in q$ then $q\cap\text{suc}(t)\ne\emptyset$;

\enskip writing $\delta_k(q) = \min\formset{{d}_t(q):t\in q\cap T_k}$
for $k\in\Bbb N$, $\lim_{k\to\infty}\delta_k(q)=\infty$.

\noindent Observe that $\delta_k(T)\ge k+1$ for each $k\in\Bbb N$, so that
$T\in\Bbb Q$ and $\Bbb Q\ne\emptyset$.

\medskip

{\bf (e)} For $q$, $q'\in\Bbb Q$ say that $q\le q'$ if $q\subseteq q'$.   Then
$(\Bbb Q,\le,T)$ is a p.o.set (that is, a pre-ordered set with a top element,
as in \KuHJ).
Observe that if $t\in q\in\Bbb Q$ then $q\cap T^{(t)}\in\Bbb Q$ and
$q\cap T^{(t)}\le q$.

\medskip

{\bf (f)} For $q$, $q'\in\Bbb Q$ and $k\in\Bbb N$ say that $q\le_k q'$
if $q\le q'$ and $q\cap T_k=q'\cap T_k$ and

\centerline{$d_{t}(q)\ge\min(k,{d}_t(q'))-2^{-k}$}

\noindent for every $t\in q$.   Note that $\le_k$ is not transitive unless
$k=0$.

\medskip

\noindent{\bf Remarks} Of course the point of the sequence $\langle n_k\rangle
_{k\in\Bbb N}$ on which the rest of this construction will depend is that it
increases `as fast as we need it to'.   The exact list given in (a) above is of
no significance and will be used only as a list of clues to the (elementary)
arguments below which depend on the rapidly-increasing nature of the
sequence.   This is why we have made no attempt to make the list as elegant
or as short as possible.

Three elements may be distinguished within the construction of $\Bbb Q$.
First, it is a p.o.set of rapidly branching trees;  that is, if $t\in q\in
\Bbb Q$, $q\cap\text{suc}(t)$ is large compared with $T_{\text{rank}
(t)}$, except for $t$ of small rank.   This is the basis of most of the
(laborious but routine) work down to 1P below.   Second, there is a natural
$\Bbb Q$-name for a subset of $X=\prod_{k\in\Bbb N}n_k$ of large measure;
a generic filter in $\Bbb Q$ leads to a branch of $T$ and hence to the
$\Psi$ of 1Q(d).   Third, the use of dp in the definition of `rapidly
branching' ((c)-(d) above) is what makes possible the side-step in the last part
of the proof of 1R.

\medskip

{\bf 1F Lemma} $\Bbb Q$ is proper.

\medskip

\noindent {\bf proof} This is a special case of Proposition 1.18 in \ShCBF.
(In fact, the arguments of 1G-1H below show that $\Bbb Q$ satisfies Axiom A, 
and is therefore proper;  see \BaHD, 2.4.)

\medskip

{\bf 1G Lemma} Let $k\in\Bbb N$ and let $\zeta$ be an ordinal.   Suppose that
$A$ is a set with $\#(A)\le\exp(2^{-k}/\gamma_j)-1$ for every $j\ge k$, and that
$\tau$ is a $\Bbb Q$-name for a member of $A$.   Let $\Delta$ be a
$\Bbb Q$-name for a countable subset of $\zeta$.   Then for every $q\in\Bbb Q$
there are a $q'\le_k q$, a function $H:T_k\to A$ and a countable (ground-model)
set $D
\subseteq\zeta$ such that

\centerline{$q'\cap T^{(t)}\Vdash_{\Bbb Q}\tau=H(t)\enskip\forall
\enskip t\in q'\cap T_k$,}

\centerline{$q'\Vdash_{\Bbb Q}\Delta\subseteq D$.}

\medskip

\noindent \bf proof (a) \rm Set $m=\#(A)$.
The point is that if $j\ge k$ and $t\in T_j$
and  $\langle S_i\rangle_{i\le m}$ is a family of subsets of $T$, then

$$\eqalign{{d}_t(\bigcup_{i\le m}S_i) &= \gamma_j\ln(\text{dp}
(\bigcup_{i\le m}\{C:t\times C\in S_i\}))\cr
&\le\gamma_j\ln(\sum_{i\le m}\text{dp}(\{C:t\times C\in S_i\}))\cr
&\le \gamma_j\ln((m+1)\max_{i\le m}\text{dp}(\{C:t\times C\in S_i\}))\cr
&=\gamma_j\ln(m+1) + \max_{i\le m}\gamma_j\ln(\text{dp}(\{C:t\times C\in S_i\}))
\cr
&\le 2^{-k}+\max_{i\le m}{d}_t(S_i).\cr}$$

\medskip

{\bf (b)} For each $a\in A$, let $S_a$ be the set

$$\eqalign{\{t:t\in q,\,&\text{rank}(t)\ge k,\,\exists\enskip p\in\Bbb Q,\,
D\in[\zeta]^{\le\omega},\cr
&p\le_k q\cap T^{(t)},\,p\Vdash_{\Bbb Q}\tau=a\enskip\&\enskip
\Delta\subseteq D\}.\cr}$$

\noindent If $t\in q\setminus S_a$ and $\text{rank}(t)\ge k$, then
$d_t(S_a)<\min(k,d_t(q))-2^{-k}$.   For if $S_a\cap\text{suc}(t)=\emptyset$,
$d_t(S_a)=-\infty$.   While if $S_a\cap\text{suc}(t)\ne\emptyset$, then for
each $s\in\text{suc}(t)\cap S_a$ we can find $p_s\in\Bbb Q$ and $D_s\in
[\zeta]^{\le\omega}$ such that $p_s\le_k q\cap T^{(s)}$, $p_s\Vdash_{\Bbb Q}
\tau=a$ and $p_s\Vdash_{\Bbb Q}\Delta\subseteq D_s$.   If we now set

\centerline{$p=\bigcup_{s\in\text{suc}(t)\cap S_a}p_s,\enskip
D=\bigcup_{s\in\text{suc}(t)\cap S_a}D_s$,}

\noindent then $p\subseteq q\cap T^{(t)}$ and $p\Vdash_{\Bbb Q}\tau=a$
and $p\Vdash_{\Bbb Q}\Delta\subseteq D$.   Because $t\notin S_a$,
$p\not\le_k q\cap T^{(t)}$ and there must be an $s\in p$
such that $d_s(p)<\min(k,d_s(q\cap T^{(t)}))-2^{-k}$;  evidently
$s=t$ and

\centerline{$d_t(S_a)<\min(k,d_t(q))-2^{-k}$,}

\noindent as claimed.

\medskip

{\bf (c)} Suppose, if possible, that there is a $t_0\in q\cap T_k\setminus
\bigcup_{a\in A}S_a$.   Set

\centerline{$p=\{t:t\in q\cap T^{(t_0)},\,t'\notin\bigcup_{a\in A}S_a
\enskip\forall\enskip t'\le t\}$.}

\noindent Then $p$ is a subtree of $T$.    For every $t\in p$ with $t\ge t_0$,

\centerline{$d_t(q)\le\max(\{d_t(p)\}\cup\{d_t(S_a):a\in A\})+2^{-k}$
}

\noindent because $\#(A)=m$.   But $d_t(S_a)<d_t(q)-2^{-k}$ for every
$a\in A$, by (b) above, so $d_t(p)\ge d_t(q)-2^{-k}$ (and $p\cap\text{suc}
(t)\ne\emptyset$).   This shows both that $p$ has no maximal elements and
that $\delta_i(p)\ge\delta_i(q)-2^{-k}$ for every $i\ge k$, so that
$p\in \Bbb Q$.   Because $\Bbb Q$ is proper, we can find a $p'\le p$ and
a countable $D\subseteq\zeta$ such that $p'\Vdash\Delta\subseteq D$ (\ShHB,
p.~81, III.1.16).   Next, there are $p''\le p'$, $a\in A$ such that
$p''\Vdash_{\Bbb Q}\tau=a$.   Let $j\in\Bbb N$ be such that $\delta_i(p'')
\ge k$ whenever $i\ge j$, and take $t\in p''$ such that $\text{rank}(t)
\ge\max(k,j)$.   Then $p''\cap T^{(t)}$ witnesses that $t\in S_a$;  which
is impossible.

\medskip

{\bf (d)} Accordingly we have for every $t\in q\cap T_k$ an $H(t)\in A$,
a countable set $D_t$ and a $p_t\in\Bbb Q$ such that

\centerline{$p_t\Vdash_{\Bbb Q}\tau=H(t)\enskip\&\enskip\Delta
\subseteq D_t$,}

\centerline{$p_t\le_k q\cap T^{(t)}$.}

\noindent Set $q'=\bigcup_{t\in q\cap T_k}p_t$, $D=\bigcup_{t\in q\cap T_k}
D_t$;  then $q'\le_k q$, $D$ is countable, $q'\Vdash_{\Bbb Q}
\Delta\subseteq D$ and $q'\cap T^{(t)}\Vdash_{\Bbb Q}
\tau=H(t)$ for every $t\in q\cap T_k$.

\medskip

{\bf 1H Lemma} Let $\langle q_k\rangle_{k\in\Bbb N}$ be a sequence
in $\Bbb Q$ such that $q_{2k+2}\le_{k+1}q_{2k+1}
\discrversionC{\break}{}\le_k q_{2k}$
for every $k\in\Bbb N$.   Then $\hat q=\bigcap_{k\in\Bbb N}q_k$ belongs
to $\Bbb Q$ and is accordingly a lower bound for $\formset{q_k:k\in\Bbb N}$
in $\Bbb Q$;  also $\hat q\cap T_{k+1}=q_{2k+1}\cap T_{k+1}$ for each $k\in
\Bbb N$.

\medskip

\noindent{\bf proof}
Because each $q_k$ is a finitely-branching subtree of $T$ with no
maximal elements, so is $\hat q$, and $d_t(\hat q)=\lim_{k\to\infty}
d_t(q_k)$ for every $t\in \hat q$.
Moreover, if $t\in\hat q$ and $k\le l\in\Bbb N$,

\centerline{$d_t(q_{2l})\ge\min(k,d_t(q_{2k}))-3.2^{-k}+3.2^{-l}$,}

\noindent (induce on $l$, using the definition of $\le_l$), so we have

\centerline{$\delta_i(\hat q)=\lim_{l\to\infty}\delta_i(q_{2l})\ge
\min(k,\delta_i(q_{2k}))-3.2^{-k}$}

\noindent for every $i$, $k\in\Bbb N$; consequently $\lim_{i\to\infty}
\delta_i(\hat q)=\infty$ and $\hat q\in\Bbb Q$.
Now if $k\in\Bbb N$ and $2k+1\le l$, $q_{l+1}\cap T_{k+1}=q_l\cap T_{k+1}$,
so $\hat q\cap T_{k+1}=q_{2k+1}\cap T_{k+1}$.

\medskip

\filename{FSp90b.tex}
\versiondate{16.9.92}

\medskip

{\bf 1I Construction: part 2}  Let $\kappa$ be the cardinal $\frak c^+$ 
(evaluated
in the ground model).

\medskip

{\bf (a)} Let $(\langle \Bbb P_{\xi}\rangle_{\xi\le\kappa},\langle\Bbb Q
_{\xi}\rangle_{\xi<\kappa})$ be a countable-support iteration of p.o.sets,
as in \KuHJ, chap. 8, such that each $\Bbb Q_{\xi}$ is a $\Bbb P_{\xi}$-name
for a p.o.set with the same definition, interpreted in
$V^{\Bbb P_{\xi}}$, as
the p.o.set $\Bbb Q$ of 1E.   (Note that $T$ is absolute, and so, in
effect, is $\langle d_t\rangle_{t\in T}$, because each $d_t$ is determined
by its values on the finite set $\Cal P(\text{suc}(t))$;  so that the
difference between $\Bbb Q$ and $\Bbb Q_{\xi}$ subsists in the power of
$\Bbb P_{\xi}$ to add new subsets of $T$.   Also each $\Bbb Q_{\xi}$
is `full' in Kunen's sense.)    
Write $\Bbb P=\Bbb P_{\kappa}$.

\medskip

{\bf (b)} If $\zeta\le\kappa$, $K\in[\zeta]^{<\omega}$, $k\in\Bbb N$ and
$p$, $p'$ belong to $\Bbb P_{\zeta}$, say that $p\le_{K,k}p'$ if
$p\le p'$ and

\centerline{$p\restr\xi\Vdash_{\Bbb P_{\xi}}p(\xi)\le_{k}p'(\xi)\enskip
\forall\enskip\xi\in K$,}

\noindent taking $\le_k$ here to be a $\Bbb P_{\xi}$-name for the relation
on $\Bbb Q_{\xi}$ corresponding to the relation $\le_k$ on $\Bbb Q$
as defined in 1E(f).   Of course $\le_{K,k}$ is not transitive unless
$K=\emptyset$ or $k=0$.

\medskip

{\bf 1J Lemma} (a) $\Bbb P_{\zeta}$ is proper for every $\zeta\le\kappa$.

(b) If $\xi<\kappa$, $\zeta\le\kappa$ then $\Bbb P_{\xi+\zeta}$ may be
identified with a dense subset of the iteration $\Bbb P_{\xi}*\Bbb P'_{\zeta}$,
where $\Bbb P'_{\zeta}$ is a $\Bbb P_{\xi}$-name with the same definition,
interpreted in $V^{\Bbb P_{\xi}}$, as the definition of $\Bbb P_{\zeta}$ in
$V$.

(c) For every $\zeta<\kappa$,

\centerline{$\Bbbone_{\Bbb P_{\zeta}}\Vdash_{\Bbb P_{\zeta}}
2^{\omega}<\kappa$.}

(d) If $\zeta\le\kappa$ has uncountable cofinality, $A$ is a (ground-model)
set, $\dot f$ is a $\Bbb P_{\zeta}$-name for a sequence in $A$ and
$p\in\Bbb P_{\zeta}$, then
we can find $\xi<\zeta$, $p'\le p$ and a $\Bbb P_{\xi}$-name $\dot g$ such that

\centerline{$p'\Vdash_{\Bbb P_{\zeta}}
\dot f=\dot g$.}

(e) If $A$ is a (ground-model) set and $\dot f$ is a $\Bbb P$-name
for a sequence in $A$, then we can find a $\xi<\kappa$ and a
$\Bbb P_{\xi}$-name $\dot g$ such that

\centerline{$\Bbbone_{\Bbb P}\Vdash_{\Bbb P}\dot f=\dot g$.}

\medskip

\noindent{\bf proof (a)} This is just because $\Bbb Q$ is proper, as noted in
1F;  see \ShHB, p.~90, Theorem III.3.2.

\medskip

{\bf (b)} This now follows by induction on $\zeta$.   The inductive step to
a successor ordinal is trivial, because if we can think of $\Bbb P_{\xi+\zeta}$
as dense in $\Bbb P_{\xi}*\Bbb P'_{\zeta}$ then we can identify $\Bbb Q_
{\xi+\zeta}$ with $\Bbb Q'_{\zeta}$.
As for the inductive step to limit $\zeta$, any member of $\Bbb P_{\xi+\zeta}$
can be regarded as $(p,p')$ where $p\in\Bbb P_{\xi}$ and 
$p'$ is a $\Bbb P_{\xi}$-name
for a member of $\Bbb P'_{\zeta}$.   On the other hand, given $(p,p')\in
\Bbb P_{\xi}*\Bbb P'_{\zeta}$, we have a $\Bbb P_{\xi}$-name $\dot J$ for
the support of $p'$ which in $V^{\Bbb P_{\xi}}$ is a countable subset of
$\zeta$.   But because $\Bbb P_{\xi}$ is proper there are a $p_1\le p$ and
a countable ground-model set $I\subseteq\zeta$ such that $p_1\Vdash
_{\Bbb P_{\xi}}\dot J\subseteq I$ (\ShHB, p.~81, III.1.16).   Now $(p_1,p')$
can be re-interpreted as a member of $\Bbb P_{\xi+\zeta}$ stronger than
$(p,p')$.   Thus $\Bbb P_{\xi+\zeta}$ is dense in $\Bbb P_{\xi}*\Bbb P'
_{\zeta}$, as claimed.

\medskip

{\bf (c)} \ShHB, p. 96, III.4.1.

\medskip

{\bf (d)} \ShHB, p. 171, V.4.4.

\medskip

{\bf (e)} By \ShHB, p. 96, III.4.1, $\Bbb P$ satisfies the
$\kappa$-c.c.;  because $\kappa$ is regular, (d) gices the result.

\medskip

{\bf 1K Definition} Let $\zeta\le\kappa$, $p\in\Bbb P_{\zeta}$.

\medskip

{\bf (a)}  Define $\bold U(p)$, $\langle p^{(\bold u)}\rangle_{\bold u\in
\bold U(p)}$ as follows.   A finite function $\bold u\subseteq\zeta\times T$
belongs to $\bold U(p)$ if

{\it either} $\bold u=\emptyset$,  in which case $p^{(\bold u)}=p$,

{\it or} $\bold u=\bold v\cup\formset{(\xi,t)}$ where $\bold v\in\bold U(p)$,
$\text{dom}(\bold v)\subseteq \xi<\zeta$, and

\centerline{$p^{(\bold v)}\restr\xi\Vdash_{\Bbb P_{\xi}}t\in p^{(\bold v)}
(\xi)$,}

\noindent in which case $p^{(\bold u)}$ is defined by writing

\centerline {$p^{(\bold u)}(\eta)=p^{(\bold v)}(\eta)\enskip\forall\enskip
\eta\in\zeta\setminus\formset{\xi}$,}

\centerline{$p^{(\bold u)}(\xi)=p^{(\bold v)}(\xi)\cap T^{(t)}$.}

\medskip

{\bf (b)} Observe that if $\bold u\in\bold U(p)$ then $p^{(\bold u)}(\xi)
=p(\xi)$ for $\xi\in\zeta\setminus\text{dom}(\bold u)$, $p^{(\bold u)}
(\xi)=p(\xi)\cap T^{(\bold u(\xi))}$
if $\xi\in\text{dom}(\bold u)$;  
$\bold U(p)$ is just the set of finite functions
$\bold u$ for which these formulae define such a $p^{(\bold u)}\in
\Bbb P_{\zeta}$.  Of course
$p^{(\bold u)}\le p$ for every $\bold u\in\bold U(p)$.

\medskip

{\bf (c)} Note that if $\xi\le\zeta$, $p\in\Bbb P_{\zeta}$, $\bold u\in
\bold U(p)$ then $\bold u\restr\xi\in\bold U(p\restr\xi)$ and
$(p\restr\xi)^{(\bold u\restr\xi)}=p^{(\bold u)}\restr\xi$.

\medskip

{\bf (d)} If $p\in\Bbb P_{\zeta}$, $\bold u\in\bold U(p)$ and
$\bold v\subseteq\zeta\times T$ is a finite function such that
$\text{dom}(\bold u)\subseteq\text{dom}(\bold v)$ and $\bold u(\xi)
\le\bold v(\xi)$ in $T$ for every $\xi\in\text{dom}(\bold u)$,
then $\bold v\in\bold U(p)$ iff $\bold v\in\bold U(p^{(\bold u)})$,
and in this case $p^{(\bold v)}=(p^{(\bold u)})^{(\bold v)}$
(induce on $\#(\bold v)$).

\medskip

{\bf (e)} We shall mostly be using not the whole of $\bold U(p)$
but the sets $\bold U(p;K,k)
\discrversionC{\break}{}=\bold U(p)\cap T_k^K$ for $K\in[\zeta]^{
<\omega}$,
$k\in\Bbb N$, writing $T_k^K$ for the set of functions from $K$ to $T_k$.

\medskip

{\bf 1L Definition} For $\zeta\le\kappa$, $K\in[\zeta]^{<\omega}$,
$k\in\Bbb N$ and $p\in\Bbb P_{\zeta}$, say that $p$ is $(K,k)$-{\bf fixed}
if for every $\eta\in K$, $\bold u\in\bold U(p;K\cap\eta,k)$ there is a
(ground-model) set $A\subseteq T_k$ such that

\centerline{$p^{(\bold u)}\restr\eta\Vdash_{\Bbb P_{\eta}}p(\eta)\cap T_k=A$.}

\noindent Equivalently, $p$ is $(K,k)$-fixed if $\bold U(p;K,k)\supseteq
\bold U(p_1;K,k)$ for every $p_1\le p$.

\medskip

{\bf 1M Lemma} 
Suppose $\zeta\le\kappa$, $K\in[\zeta]^{<\omega}$, $k\ge 1$ and that $A$
is a finite set with $2^{c^m}a^{c^{m-1}}\le\exp(2^{-k}/\gamma_i)$ for every
$i\ge k$, where $c=\#(T_k)$, $m=\#(K)$ and $a=\#(A)$.   Let $\tau$ be a 
$\Bbb P_{\zeta}$-name for a member of $A$, and $\Delta$ a
$\Bbb P_{\zeta}$-name for a countable subset of $\kappa$.   Then for any
$p\in\Bbb P_{\zeta}$ there are $p_1\le_{K,k}p$, a function $H:T_k^K\to A$
and a countable (ground-model) set $D\subseteq\kappa$ such that

\centerline{$p_1$ is $(K,k)$-fixed,}

\centerline{$p_1^{(\tbf{u})}
\Vdash_{\Bbb P_{\zeta}}\tau=H(\tbf{u})\enskip\forall\enskip\tbf{u}\in\tbf{U}(p_1;K,k)$,
}

\centerline{$p_1\Vdash_{\Bbb P_{\zeta}}\Delta\subseteq D$.}

\medskip

\noindent{\bf proof} Induce on $m=\#(K)$.   If $m=0$ we may take any
$a\in A$, $p_1'\le p$ such that $p_1'\Vdash\tau=a$, and (again using \ShHB,
III.1.16, this time based on 1Ja) 
a countable $D$ and a $p_1\le p_1'$ such that $p_1\Vdash_{\Bbb P_
{\zeta}}\Delta\subseteq D$;  now set $H(\emptyset)=a$.

For the inductive step to $\#(K)=m\ge 1$, let $\xi$ be $\max K$.   As explained
in 1Jb, $\Bbb P_{\zeta}$ may be regarded as a dense subset of $\Bbb P_{\xi+1}
*\Bbb P'$;  arguing momentarily in $V^{\Bbb P_{\xi+1}}$ we can find
a $\Bbb P_{\xi+1}$-name $\hat r_0$ for a member of
$\Bbb P'$, a $\Bbb P_{\xi+1}$-name $\tau'$ for a member of $A$ and a
$\Bbb P_{\xi+1}$-name  $\Delta'$ for  a countable set such that

\centerline{$(p\restr\xi+1,\hat r_0)\le p$ in $\Bbb P_{\xi+1}*\Bbb P'$,}

\centerline{$(p\restr\xi+1,\hat r_0)\Vdash_{\Bbb P_{\xi+1}*\Bbb P'}\tau'=\tau$,}

\centerline{$(p\restr\xi+1,\hat r_0)\Vdash_{\Bbb P_{\xi+1}*\Bbb P'}
\Delta\subseteq\Delta'$.}

\noindent Now let $\Delta_0'$ be a $\Bbb P_{\xi+1}$-name for a countable
subset of $\zeta\setminus(\xi+1)$ such that

\centerline{$\Bbbone_{\Bbb P_{\xi+1}}\Vdash_{\Bbb P_{\xi+1}}
\text{supp}(\hat r_0)=\Delta_0'$.}

\noindent Because $\#(A)=a<2^{c^m}a^{c^{m-1}}\le\exp(2^{-k}/\gamma_i)$ for
every $i\ge k$, we can use Lemma 1G in $V^{\Bbb P_{\xi}}$ to find
$\tilde H$, $\tilde q$, $\tilde\Delta$ such that

\centerline{$\tilde H$ is a $\Bbb P_{\xi}$-name for a function from
$T_k$ to $A$,}

\centerline{$\tilde\Delta$ is a $\Bbb P_{\xi}$-name for a countable subset
of $\kappa$,}

\centerline{$\tilde q\in\Bbb Q_{\xi}$,}

\centerline{$p\restr\xi\Vdash_{\Bbb P_{\xi}}\tilde q\le_k p(\xi)$,}

\centerline{$p\restr\xi\Vdash_{\Bbb P_{\xi}}\bigl(\tilde q\cap T^{(t)}\Vdash
_{\Bbb Q_{\xi}}\tau'=\tilde H(t)\enskip\forall\enskip t\in\tilde q\cap
T_k\bigr)$,}

\centerline{$p\restr\xi\Vdash_{\Bbb P_{\xi}}\bigl(\tilde q\Vdash_{\Bbb Q_{\xi}}
\Delta'\cup\Delta'_0\subseteq\tilde\Delta\bigr)$.}

Now consider the pair $(\tilde H,\tilde q\cap T_k)$.   This can be regarded
as a $\Bbb P_{\xi}$-name for a member of $A_1=A^{T_k}\times\Cal PT_k$, and
$a_1=\#(A_1)=2^ca^c$, so

\centerline{$2^{c^{m-1}}a_1^{c^{m-2}}=2^{2c^{m-1}}a^{c^{m-1}}
\le 2^{c^m}a^{c^{m-1}}\le \exp(2^{-k}/\gamma_j)\enskip\forall
\enskip i\ge k$.}

\noindent The inductive hypothesis therefore tells us that there are
$\hat p_1\le_{K\cap\xi,k}p\restr\xi$, $H^*:T_k^{K\cap\xi}\to A^{T_k}$,
$F^*:T_k^{K\cap\xi}\to\Cal PT_k$ and a countable $D\subseteq\kappa$ such
that

\centerline{$\hat p_1$ is $(K\cap\xi,k)$-fixed,}

\centerline{$\hat p_1^{(\tbf{u})}\Vdash_{\Bbb P_{\xi}}\tilde H=H^*(\tbf{u})
\enskip\&\enskip\tilde q\cap T_k=F^*(\tbf{u})$}

\noindent for every $\tbf{u}\in\tbf{U}(\hat p_1;K\cap\xi,k)$, and

\centerline{$\hat p_1\Vdash_{\Bbb P_{\xi}}\tilde\Delta\subseteq D$.
}

At this point we observe that 

\centerline{$\hat p_1\Vdash_{\Bbb P_{\xi}}\bigl(\tilde q\Vdash_{\Bbb Q_{\xi}}
\text{supp}(\hat r_0)\subseteq D\bigr)$.}

\noindent Now the only difference between $\Bbb P_{\xi+1}*\Bbb P'$ and
$\Bbb P_{\zeta}$ is that for members of the former their supports have to
be regarded as $\Bbb P_{\xi+1}$-names for countable subsets of $\zeta$,
and these are not always reducible to countable ground-model sets.   But in
the present case this difficulty does not arise and we have a $p_1\in
\Bbb P_{\zeta}$ defined by saying that $p_1\restr\xi=\hat p_1$,
$p_1\restr\xi\Vdash_{\Bbb P_{\xi}}p_1(\xi)=\tilde q$, and $p_1\restr\eta
\Vdash_{\Bbb P_{\eta}}p_1(\eta)=\hat r_0(\eta)$ for $\xi<\eta<\zeta$;  then
$\text{supp}(p_1)\subseteq\text{supp}(\hat p_1)\cup\{\xi\}\cup(D\cap\zeta)$
is countable.

Consequently $p_1\in\Bbb P_{\zeta}$ is well-defined and now, setting
$H(\tbf{u}^{\smallfrown}t)=H^*(\tbf{u})(t)$ for $\tbf{u}\in T_k^{K\cap\xi}$,
$t\in T_k$,

\centerline{$p_1\le_{K,k}p$,}

\centerline{$p_1\Vdash_{\Bbb P}\Delta\subseteq D$,}

\centerline{$\tbf{U}(p_1;K,k)=\{\tbf{u}^{\smallfrown}t:
\tbf{u}\in\tbf{U}(\hat p_1;K\cap\xi,k),\,t\in F^*(\tbf{u})\}$,}

\centerline{$p_1^{(\tbf{v})}\Vdash_{\Bbb P_{\zeta}}\tau=
H(\tbf{v})\enskip\forall\enskip\tbf{v}\in\tbf{U}(p_1;K,k)$}

\noindent and finally

\centerline{$p_1^{(\tbf{u})}\restr\xi\Vdash_{\Bbb P_{\xi}}
p_1(\xi)\cap T_k=F^*(\tbf{u})\enskip\forall\enskip\tbf{u}\in\tbf{U}(p_1;
K\cap\xi,k)$,}

\noindent so that $p_1$ is $(K,k)$-fixed, and the induction proceeds.

\medskip

{\bf 1N Lemma} Suppose $\zeta\le\kappa$, $\langle K_k\rangle_{k\in\Bbb N}$
is an increasing sequence of finite subsets of $\zeta$, $\langle p_k\rangle
_{k\in\Bbb N}$ is a sequence in $\Bbb P_{\zeta}$;  suppose that

\centerline{$p_{2k+2}\le_{K_k,k+1}p_{2k+1}\le_{K_k,k}p_{2k}$}

\noindent for every $k\in\Bbb N$ and that $\bigcup_{k\in\Bbb N}\text{supp}
(p_k)\subseteq\bigcup_{k\in\Bbb N}K_k$.   Then there is a $\hat p\in\Bbb P_
{\zeta}$ such that $\hat p\le p_k$ for every $k\in\Bbb N$,
$\text{supp}(\hat p)\subseteq\bigcup_{k\in\Bbb N}K_k$ and

\centerline{$\hat p\restr\xi\Vdash_{\Bbb P_{\xi}}\hat p(\xi)\cap T_k
=p_{2k+1}\cap T_k\enskip\forall\enskip k\in\Bbb N,\,\xi\in K_k$,}

\centerline{$\hat p\restr\xi\Vdash_{\Bbb P_{\xi}}\hat p(\xi)\cap T_{k+1}
=p_{2k+2}\cap T_{k+1}\enskip\forall\enskip k\in\Bbb N,\,\xi\in K_{k+1}$,}

\noindent so that

\centerline{$\bold U(\hat p;K_k,k)\supseteq\bold U(p_{2k+1};K_k,k)$ and
$\bold U(\hat p;K_k,k+1)\supseteq\bold U(\hat p_{k+2};K_k,k+1)$}

\noindent for every $k\in\Bbb N$.

\medskip

\noindent{\bf proof} For each $\xi<\zeta$ choose $\hat p(\xi)$ such that

\centerline{$\Bbbone_{\Bbb P_{\xi}}\Vdash_{\Bbb P_{\xi}}\hat p(\xi)
=\bigcap_{k\in\Bbb N}p_k(\xi)$.}

\noindent An easy induction on $\xi$ shows that 
$\hat p\restr\xi\in\Bbb P_{\xi}$
for every $\xi\le\zeta$; for if $\xi\in\zeta\setminus\bigcup_{k\in\Bbb N}K_k$
then 

\centerline{$\Bbbone_{\Bbb P_{\xi}}\Vdash_{\Bbb P_{\xi}}\hat p(\xi)=T
=\Bbbone_{\Bbb Q_{\xi}}$,}

\noindent while if $k\in\Bbb N$ and $\xi\in K_k$ then

\centerline{$\hat p\restr\xi\Vdash_{\Bbb P_{\xi}}p_{2l+2}(\xi)\le_{l+1}
p_{2l+1}(\xi)\le_l p_{2l}(\xi)\enskip\forall\enskip l\ge k$,}

\noindent so that by Lemma 1H,

$$\eqalign{\hat p\restr\xi\Vdash_{\Bbb P_{\xi}}\hat p(\xi)\in\Bbb Q_{\xi}
\,&\&\,\hat p(\xi)\cap T_{k+1}=p_{2k+1}\cap T_{k+1}=p_{2k+2}\cap T_{k+1}\cr
&\&\,\hat p(\xi)\cap T_k=p_{2k+1}\cap T_k.\cr}$$

\noindent It follows at once that $\bold U(\hat p;K_k,k)\supseteq\bold U(
p_{2k+1};K_k,k)$, $\bold U(\hat p;K_k,k+1)\supseteq
\bold U(p_{2k+2};K_k,k+1)$ for
every $k\in\Bbb N$.

\medskip

{\bf 1O Lemma}
Suppose that $0<\zeta\le\kappa$, $\sigma$ is a $\Bbb P_{\zeta}$-name for
a member of $\prod_{k\in\Bbb N}n_k$, and $p\in\Bbb P_{\zeta}$.   Then
we can find a $\hat p$ and sequences  $\langle K_k\rangle_{k\in\Bbb N}$,
$\langle H_k\rangle_{k\in\Bbb N}$ such that

\quad $\hat p\in\Bbb P_{\zeta}$, $\hat p\le p$;

\quad $\langle K_k\rangle_{k\in\Bbb N}$ is an increasing sequence of subsets
of $\zeta$, $\#(K_k)\le k+1$ for every $k$, $K_0=\formset{0}$;

\quad$\text{supp}(\hat p)\subseteq\bigcup_{k\in\Bbb N}K_k$;

\quad $\hat p$ is $(K_k,k)$-fixed and $(K_k,k+1)$-fixed for every $k$;

\quad $H_k$ is a function from $T_{k+1}^{K_k}$ to $n_k$ for every $k$;

\quad $\hat p^{(\bold u)}\Vdash_{\Bbb P_{\zeta}}\sigma(k)=H_k(\bold u)$
whenever $k\in\Bbb N$ and $\bold u\in\bold U(\hat p;K_k,k+1)$.

\medskip

\noindent{\bf proof} Using Lemma 1M, we can find sequences $\langle p_k\rangle
_{k\in\Bbb N}$, $\langle K_k\rangle_{k\in\Bbb N}$ and 
\discrversionC{\break}{}$\langle H_k\rangle
_{k\in\Bbb N}$ such that

\quad $p=p_0$, $K_0=\formset{0}$;

\quad $\#(K_{k+1})\le k+2$, $K_{k+1}\supseteq K_k$;

\quad $p_{2k+1}\le_{K_k,k}p_{2k}$, $p_{2k+1}$ is $(K_k,k)$-fixed;

\quad $H_k:T_{k+1}^{K_k}\to n_k$ is a function;

\quad $p_{2k+2}\le_{K_k,k+1}p_{2k+1}$, $p_{2k+2}$ is $(K_k,k+1)$-fixed,

\qquad $p_{2k+2}^{(\bold u)}\Vdash_{\Bbb P_{\zeta}}\sigma(k)=H_k(\bold u)
\enskip\forall\enskip\bold u\in\bold U(p_{2k+2};K_k,k+1)$

\noindent for every $k\in\Bbb N$.   Furthermore, we may do this in such a way
that 
$\bigcup_{k\in\Bbb N}K_k$ includes
$\bigcup_{k\in\Bbb N}\text{supp}(p_k)$.
We need of course to know that the $n_k$ are rapidly increasing; specifically,
that

\centerline{$2^{c_k^{k+1}}\le\exp(2^{-k}/\gamma_i)\enskip\forall\enskip i\ge k$}

\noindent (when choosing $p_{2k+1}$) and that

\centerline{$
2^{c_{k+1}^{k+1}}n_k^{c_{k+1}^k}\le\exp(2^{-k-1}/\gamma_i)\enskip\forall
\enskip i\ge k+1$}

\noindent (when choosing $p_{2k+2}$), where we write $c_k=\#(T_k)$.
But as $c_k
\le\prod_{i<k}2^{n_i}$,
this is a consequence of 1E(a)(i) and (iii).

Armed with the sequences $\langle p_k\rangle_{k\in\Bbb N}$, $\langle K_k\rangle
_{k\in\Bbb N}$ we may now use Lemma 1N to find a $\hat p$ as described there.
Because $p_{2k+1}$ is $(K_k,k)$-fixed and
\discrversionC{\break}{}
$\bold U(\hat p;K_k,k)\supseteq\bold U(p_{2k+1};K_k,k)$ we must have equality
here and
$\hat p$ is $(K_k,k)$-fixed for every $k\in\Bbb N$.
Similarly, $\hat p$ is $(K_k,k+1)$-fixed for every $k$.   Moreover,
if $\bold u\in\bold U(\hat p;K_k,k+1)=\bold U(p_{2k+2};K_k,k+1)$ we have
$\hat p^{(\bold u)}\le p_{2k+2}^{(\bold u)}$, so

\centerline{$\hat p^{(\bold u)}\Vdash_{\Bbb P_{\zeta}}\sigma(k)=H_k(\bold u)$}

\noindent  as
required.

\medskip

{\bf 1P Lemma} Suppose that $\zeta\le\kappa$, $p\in\Bbb P_{\zeta}$,
$k\in\Bbb N$, $K\in[\zeta]^{<\omega}$ and $\bold V$ is a non-empty subset
of $\bold U(p;K,k)$.   Then we have a $p_1=\bigvee_{\bold v\in\bold V}
p^{(\bold v)}$ defined (up to $\le$-equivalence in $\Bbb P_{\zeta}$) by
saying

\quad if $\xi\in\zeta\setminus K$ then $p_1(\xi)=p(\xi)$;

\quad if $\xi\in K$ then

\centerline{$(p_1\restr\xi)^{(\bold u)}\Vdash_{\Bbb P_{\xi}}
p_1(\xi)=\bigcup\formset{p(\xi)\cap T^{(t)}:\exists\enskip\bold v\in\bold V$
such that $\bold v\restr\xi+1=\bold u^{\smallfrown}t}$}

\quad for $\bold u\in\formset{\bold v\restr\xi:\bold v\in\bold V}$.

\noindent Now $p_1\le p$ and 
if $\xi<\zeta$, $t\in p_1(\xi)$, $\text{rank}(t)\ge k$ we shall have

\centerline{$p_1\restr\xi\Vdash_{\Bbb P_{\xi}}\text{suc}(t)\cap p_1(\xi)
=\text{suc}(t)\cap p(\xi)$;}

\noindent so if $\xi<\zeta$, $i\ge k$ we have

\centerline{$p_1\restr\xi\Vdash_{\Bbb P_{\xi}}\delta_i(p_1(\xi))\ge
\delta_i(p(\xi))$.}

\noindent If $p_2\le p_1$ there is some $\bold v\in
\bold V$ such that $p_2$ is compatible with $p_1^{(\bold v)}=p^{(\bold v)}$.
If $k\le l\in\Bbb N$, $K\subseteq L\in[\zeta]^{<\omega}$ then

\centerline{$\bold U(p_1;L,l)=\formset{\bold w:\bold w\in\bold U(p;L,l),
\,\exists\enskip \bold v\in\bold V$ such that 
$\bold v(\xi)\le\bold w(\xi)\enskip\forall\enskip\xi\in K}$,}

\noindent and $p_1^{(\bold w)}=p^{(\bold w)}$ for every $\bold w\in
\bold U(p_1;L,l)$; consequently, $p_1$ is $(L,l)$-fixed if $p$ is.

\medskip

\noindent {\bf proof} Requires only a careful reading of the definitions.

\medskip

\noindent{\bf Remark} Note that 1G-1P are based just on the fact that $\Bbb Q$
is a p.o.set of rapidly branching trees;  the exact definition of `rapidly
branching' in 1E(c) is relevant only to some of the detailed calculations.
Similar ideas may be found in \BJSpHI and \ShCBF.

\medskip

{\bf 1Q Construction:  part 3 (a)} Set $X=\prod_{k\in\Bbb N}n_k$.   Then
$X$, with its product topology, is a compact metric space.   Let $\mu$ be
the natural Radon probability on $X$, the product of the uniform probabilities
on the factors.

\medskip

{\bf (b)} For each $k\in\Bbb N$ set $l_k=\lceil(\ln n_k)^2\rceil$.   Take
$W'_k\subseteq n_k\times n_k$ such that $\#(W'_k)\le 2^{-k-1}n_k^2$ and
whenever $I\in[n_k]^{l_k}$ and $J_0,\ldots,J_{l_k-1}$ are disjoint members of
$[n_k]^{\le k}$, there are $i\in I$ and $j<l_k$ such that $\formset{i}
\times J_j\subseteq W'_k$.   (This is possible by Lemma 1B and 1E(a)(iv).)
Set $W_k=W'_k\cup\formset{(i,i):i<n_k}$.

Write $R$ for

\centerline{$\formset{(x,y):x,\,y\in X,\,(x(k),y(k))\in W_k\enskip\forall
\enskip k\in\Bbb N,\,\formset{k:x(k)=y(k)}\text{ is finite}}$;}

\noindent then $R$ is negligible for
the product measure of  $X\times X$.  
For $r\in\Bbb N$ write $R_r$ for the set

\centerline{$\formset{(x,\langle y_i\rangle_{i<r}):x\in X,\,(x,y_i)\in R\enskip
\forall\enskip i<r}\subseteq X\times X^r$.}

\noindent We shall frequently wish to interpret the formulae for the
sets $X$, $R_r$ in
$V^{\Bbb P}$;  when doing so we will write
$\ulcorner X\urcorner$, $\ulcorner R_r\urcorner$.

\medskip

{\bf (c)}   
Write

\centerline
{$\Cal L=\formset{\langle L_k\rangle_{k\in\Bbb N}:L_k\subseteq n_k\enskip
\forall\enskip k\in\Bbb N,\,\prod_{k\in\Bbb N}\#(L_k)/n_k>0}$.}

\noindent Again, we shall wish to distinguish between the ground-model
set $\Cal L$ and a corresponding $\Bbb P$-name $\ulcorner \Cal L\urcorner$.

\medskip

{\bf (d)} For each $k\in\Bbb N$ let $\Phi_k$ be the $\Bbb P$-name 
for a subset of $n_k$ defined (up to equivalence) by saying that

\centerline{$p\Vdash_{\Bbb P}\Phi_k=C_k(t)$ }

\noindent whenever $\text{rank}(t)> k$ and
$p(0)\subseteq T^{(t)}$.
(Here $C_k(t)$ is the $k$th factor of $t$, as described in 1E(b).)
Let $\Psi_k$, $\Psi$ be $\Bbb P$-names for the subsets of
$\ulcorner X\urcorner$ given by

\centerline{$\Bbbone_{\Bbb P}\Vdash_{\Bbb P}\Psi_k=\formset{\sigma:
\sigma\in\ulcorner X\urcorner,\,
\sigma(i)\in\Phi_i\enskip\forall\enskip i\ge k},\,\Psi=\bigcup_{k\in\Bbb N}
\Psi_k$.}

\noindent Then we have

\centerline{$\Bbbone_{\Bbb P}\Vdash_{\Bbb P}\#(\Phi_k)\ge(1-2^{-k-1})n_k\enskip
\forall\enskip k\in\Bbb N$,}

\noindent so that

\centerline{$\Bbbone_{\Bbb P}\Vdash_{\Bbb P}\ulcorner\mu\urcorner(\Psi)=1$.}

\medskip

{\bf 1R Main Lemma} If $r\in\Bbb N$ and $D\subseteq X^{r}$ is a 
(ground-model)   
set such that $D\cap(\prod_{k\in\Bbb N}L_k)^r\ne\emptyset$ for every 
(ground-model)
sequence $\langle L_k\rangle_{k\in\Bbb N}\in\Cal L$, then for every
(ground-model) sequence $\langle L_k\rangle_{k\in\Bbb N}\in\Cal L$

\centerline{$\Bbbone_{\Bbb P}\Vdash_{\Bbb P}\Psi\cap\prod_{k\in\Bbb N}L_k
\subseteq \hbox{$\ulcorner R_r\urcorner$}^{-1}
[D\cap(\prod_{k\in\Bbb N}L_k)^r]$.}

\medskip

\noindent{\bf proof (a)}   
Let $\langle L_k\rangle_{k\in\Bbb N}\in\Cal L$, let
$\sigma$ be a $\Bbb P$-name such that

\centerline{$\Bbbone_{\Bbb P}\Vdash_{\Bbb P}\sigma\in\Psi\cap\prod_{k\in\Bbb N}
L_k$,}

\noindent and let $p\in\Bbb P$.
Write $D'$ for $D\cap(\prod_{k\in\Bbb N}L_k)^r$.
Let $k_0\ge r$, $p_1\le p$ be such that
$p_1\Vdash_{\Bbb P}\sigma\in\Psi_{k_0}$.   By Lemma 1O,
we have a $p_2\le p_1$, an increasing sequence 
$\langle K_k\rangle_{k\in\Bbb N}$
of finite subsets of $\kappa$, and a sequence $\langle H_k\rangle_{k\in\Bbb N}$
of functions such that

\quad$p_2$ is $(K_k,k)$-fixed and $(K_k,k+1)$-fixed for every $k\in\Bbb N$,

\quad$p_2^{(\bold u)}\Vdash_{\Bbb P}\sigma(k)=H_k(\bold u)$
whenever $\bold u\in\bold U(p_2;K_k,k+1),\,k\in\Bbb N$;

\quad$\bigcup_{k\in\Bbb N}K_k\supseteq\text{supp}(p_2)$;

\quad$0\in K_0$, $\#(K_k)\le k+1$ for every $k\in\Bbb N$.

\medskip

{\bf (b)} For $k\ge k_0$, 
let $Z_k$ be the cartesian product set $n_k^{r}$
and take $\tilde W_k$ to be

\centerline{$\formset{(i,z):i<n_k,\,z\in Z_k,\,(i,z(j))\in W_k\enskip\forall
\enskip j<r}$.}

\noindent Set $\Cal A_k=\Cal Pn_k\setminus\formset{\emptyset}$ and
$\Bbb T_k=\Cal P(\Cal A_k^{K_k})\setminus\formset{\emptyset}$;  define
$\preccurlyeq_k$ on $\Bbb T_k$ as in Lemma 1C, taking $l_k$ and $K_k$
(with the order induced by that of $\kappa$) in
place of $l$ and $m$ there.

For each $\bold u\in\bold U(p_2;K_k,k)$ set

\centerline{$\Cal T_{\bold u}=\formset{\bold c:\bold u^{\wedge}\bold c\in
\bold U(p_2;K_k,k+1)}\in\Bbb T_k$,}

\noindent where for $\bold u\in T_k^K$, $\bold c\in(\Cal Pn_k)^K$ we write

\centerline{$\bold u^{\wedge}\bold c=\langle\bold u(\xi)\times\bold c
(\xi)\rangle_{\xi\in K}$.}

\noindent By Corollary 1D, 
we may find for each such $\bold u$ a $w_{\bold u}<n_k$
and a set $J_{\bold u}\subseteq n_k$ such that $\#(J_{\bold u})\le rl_k$ and

{\it either} there is a $\Cal T\preccurlyeq_k \Cal T_{\bold u}$ such that
$H_k(\bold u^{\wedge}\bold c)=w_{\bold u}$ for every $\bold c\in\Cal T$

{\it or} for every $z\in(n_k\setminus J_{\bold u})^r$ there is a $\Cal T
\preccurlyeq_k\Cal T_{\bold u}$ such that $(H_k(\bold u^{\wedge}\bold c),
z)\in\tilde W_k$, that is, 

\centerline{$(H_k(\bold u^{\wedge}\bold c),z(j))\in W_k
\enskip\forall\enskip j<r$,}

\noindent for every $\bold c\in\Cal T$.

Set 

\centerline{$\tilde I_k=\formset{w_{\bold u}:\bold u\in\bold U(p_2;K_k,k)}$,}

\centerline{$\tilde J_k=\bigcup\formset{J_{\bold u}:\bold u\in
\bold U(p_2;K_k,k)},$}

\noindent so that

\centerline{$\#(\tilde I_k)\le
\#(\bold U(p_2;K_k,k))\le\#(T_k^{K_k})\le(\prod_{i<k}2^{n_i})^{k+1}$,}

\centerline{$\#(\tilde J_k)\le kl_k(\prod_{i<k}2^{n_i})^{k+1}\le 2^{-k}n_k$}

\noindent by 1E(a)(vi).

\medskip

{\bf (c)} Let $k_1\ge\max(k_0,1)$ be such that

\centerline{$\delta_k(p_2(0))\ge 2$, $\tilde J_k\not\supseteq L_k$}

\noindent for every $k\ge k_1$.   Take any $\bold v^*\in\bold U(p_2;K_{k_1-1},
k_1)$ and set $p_3=p_2^{(\bold v^*)}$.
Then $p_3$ is $(K_k,k)$-fixed and $(K_k,k+1)$-fixed for every $k\ge k_1$,
and

\centerline{$\Cal T_{\bold u}=\formset{\bold c:\bold u^{\wedge}\bold c\in
\bold U(p_3;K_k,k+1)}$}

\noindent whenever $k\ge k_1$, $\bold u\in\bold U(p_3;K_k,k)$.

\medskip

{\bf (d)} We have

\centerline{$p_3\Vdash_{\Bbb P}\sigma(i)=H_i(\bold v_i^*)\enskip\forall
\enskip i<k_1$,}

\noindent where $\bold v_i^*$ is that member of $T_{i+1}^{K_i}$ such that
$\bold v_i^*(\eta)\le\bold v^*(\eta)$ for every $\eta\in K_i$.   Set
$L'_k=\formset{H_k(\bold v_k^*)}$ for $k<k_1$, $L'_k=L_k\setminus \tilde
J_k$ for $k\ge k_1$;  then $\prod_{k\in\Bbb N}\#(L'_k)/n_k>0$, because
$\prod_{k\in\Bbb N}\#(L_k)/n_k>0$ and $\sum_{k\in\Bbb N}\#(\tilde J_k)
/n_k<\infty$.
So there is a $\tilde z\in D\cap(\prod_{k\in\Bbb N}L'_k)^r\subseteq D'$.   
Writing
$z_k=\langle\tilde z(j)(k)\rangle_{j<r}$ for $k\in\Bbb N$, we have
$z_k\in (n_k\setminus \tilde J_k)^r$ for $k\ge k_1$, and

\centerline{$p_3\Vdash_{\Bbb P}(\sigma(k),z_k)\in\tilde W_k$}

\noindent for $k<k_1$, because $(i,i)\in W_k$ for $i<n_k$.

\medskip

{\bf (e)} For each $k\ge k_1$, $\bold u\in\bold U(p_3;K_k,k)$ choose
$\Cal T'_{\bold u}
\preccurlyeq_k\Cal T_{\bold u}\in\Bbb T_k$ such that

{\it either} $H_k(\bold u^{\wedge}\bold c)=w_{\bold u}\in\tilde I_k$ for
every $\bold c\in\Cal T'_{\bold u}$

{\it or} $(H_k(\bold u^{\wedge}\bold c),z_k)\in\tilde W_k$ for every
$\bold c\in \Cal T'_{\bold u}$.

\noindent Define $\langle \Cal S_k\rangle_{k\ge k_1}$, $\langle\tilde p_k
\rangle_{k\ge k_1}$ by

\quad $\tilde p_{k_1}=p_3$,

\quad $\Cal S_k=\formset{\bold u^{\wedge}\bold c:\bold u\in\bold U(\tilde
p_k;K_k,k)$, $\bold c\in\Cal T'_{\bold u}}$,

\quad $\tilde p_{k+1}=\bigvee\formset{\tilde p_k^{(\bold v)}:\bold v\in
\Cal S_k}$

\noindent for every $k\ge k_1$, as in Lemma 1P.   An easy induction on $k$
shows that

\centerline{$\Cal T'_{\bold u}=\formset{\bold c:\bold u^{\wedge}\bold c
\in\bold U(\tilde p_{k+1};K_k,k+1)}$}

\noindent
whenever
$\bold u\in\bold U(\tilde p_k;K_k,k)$, $k\ge k_1$,
that $\tilde p_k$ is $(K_l,l)$-fixed and $(K_l,l+1)$-fixed whenever
$k_1\le k\le l$, and that $\tilde p_k^{(\bold v)}=p_3^{(\bold v)}$
whenever $k_1\le k\le l$ and $\bold v\in\bold U(\tilde p_k;K_l,l)\cup
\bold U(\tilde p_k;K_l,l+1)$.
Also $\text{supp}(\tilde p_k)\subseteq\bigcup_{l\in\Bbb N}K_l$ for
every $k\ge k_1$.

\medskip

{\bf (f)} It is likewise easy to see that, for $k\ge k_1$,

\quad$\tilde p_{k+1}\le\tilde p_k$,

\quad$\tilde p_{k+1}\restr\xi\Vdash_{\Bbb P_{\xi}}\tilde p_{k+1}(\xi)\cap
T_k=\tilde p_k(\xi)\cap T_k\enskip\forall\enskip\xi<\kappa$,

\quad$\tilde p_{k+1}\restr\xi\Vdash_{\Bbb P_{\xi}}
\tilde p_{k+1}(\xi)\cap\text{suc}(t)=\tilde p_k(\xi)\cap\text{suc}(t)\enskip
\forall\enskip t\in\tilde p_{k+1}(\xi)\cap T_i$

\noindent unless $i=k$ and $\xi\in K_k$,

\quad$\tilde p_{k+1}\Vdash_{\Bbb P}\sigma(k)\in\tilde I_k$ {\it or} 
$(\sigma(k),
z_k)\in\tilde W_k$.

\medskip

{\bf (g)} On the other hand, if $k\ge k_1$ and $\xi\in K_k$,

$$\eqalign{\tilde p_{k+1}\restr\xi\Vdash_{\Bbb P_{\xi}}
\text{dp}(\formset{C:t\times C\in\tilde p_{k+1}(\xi)})&\ge
\text{dp}(\formset{C:t\times C\in\tilde p_k(\xi)})/2l_k\cr
&\forall\enskip t\in\tilde p_{k+1}(\xi)\cap T_k.\cr}$$

\noindent To see this, take any $q\le\tilde p_{k+1}\restr\xi$ and $t$
such that 

\centerline{$q\Vdash_{\Bbb P_{\xi}}t\in\tilde p_{k+1}(\xi)\cap T_k
=\tilde p_k(\xi)\cap T_k$.}

\noindent We may suppose that $\bold v_0\in\Cal S_k$ is such that
$q\le\tilde p_k^{(\bold v_0)}\restr\xi=\tilde p_{k+1}^{(\bold v_0)}\restr\xi
=q_1$.   Now $\tilde p_{k+1}$ is $(K_k,k+1)$-fixed so there must be a
$t'\ge t$ such that

\centerline{$q_1\Vdash_{\Bbb P_{\xi}}t'\in T_{k+1}\cap\tilde p_{k+1}(\xi)$.}

\noindent There is accordingly a $\bold v_1\in\Cal S_k$ such that
$\bold v_0\restr\xi=\bold v_1\restr\xi$ and $\bold v_1(\xi)=t'$.
Express $\bold v_1$ as $\bold u^{\wedge}\bold c_1$ where $\bold u\in
\bold U(\tilde p_k;K_k,k)$ and $\bold c_1\in\Cal T'_{\bold u}$.
Of course $\bold u(\xi)=t$.

Now

\quad$q_1\Vdash_{\Bbb P_{\xi}}\formset{C:t\times C\in\tilde p_{k+1}(\xi)}
\supseteq\formset{\bold c(\xi):\bold c\in\Cal T'_{\bold u},\,
\bold c\restr\xi=\bold c_1\restr\xi}$,

\quad$q_1\Vdash_{\Bbb P_{\xi}}\formset{C:t\times C\in\tilde p_k(\xi)}
=\formset{\bold c(\xi):\bold c\in\Cal T_{\bold u},\,\bold c\restr\xi
=\bold c_1\restr\xi}$

\noindent because $\tilde p_k$ and $\tilde p_{k+1}$ are both $(K_k,k+1)$-fixed,
while $\bold v_0\restr\xi=(\bold u\restr\xi)^{\wedge}(\bold c_1\restr\xi)$.
But because $\Cal T'_{\bold u}\preccurlyeq_k\Cal T_{\bold u}$,

\centerline{$\text{dp}(\formset{\bold c(\xi):\bold c\in\Cal T'_{\bold u},
\,\bold c\restr\xi=\bold c_1\restr\xi})\ge\text{dp}(\formset{\bold c(\xi):
\bold c\in\Cal T_{\bold u},\,\bold c\restr\xi=\bold c_1\restr\xi})/2l_k$.}

\noindent So we get

\centerline{$q\le q_1\Vdash_{\Bbb P_{\xi}}
\text{dp}(\formset{C:t\times C\in\tilde p_{k+1}(\xi)})\ge
\text{dp}(\formset{C:t\times C\in\tilde p_k(\xi)})/2l_k$.}

\noindent As $q$ and $t$ are arbitrary, we have the result.

\medskip

{\bf (h)} Because $\gamma_k\ln(2l_k)\le 2^{-k}$ (by 1E(a)(v)),

\centerline{$\tilde p_{k+1}\le_{K_k,k}\tilde p_k$}

\noindent for every $k\ge k_1$.
Also, supp$(\tilde p_k)\subseteq\bigcup_{l\in\Bbb N}K_l$ for every
$k\ge k_1$.  By Lemma 1N, there is a $p_4\in\Bbb P$ such that
$p_4\le\tilde p_k$ for every $k\ge k_1$.
Moreover, we may take it that

\centerline{$p_4(0)=\bigcap_{k\ge k_1}\tilde p_k(0)$}

\noindent (as in Lemma 1H), so that

\centerline{$\delta_i(p_4(0))\ge\delta_i(p_3(0))-\gamma_i\ln(2l_i)\ge 1$}

\noindent whenever $i\ge k_1$.   Note that for $k\ge k_1$,

\centerline{$p_4\Vdash_{\Bbb P}\sigma(k)\in\tilde I_k$ {\it or}
$(\sigma(k),z_k)\in\tilde W_k$,}

\noindent while for $k<k_1$,

\centerline{$p_4\Vdash_{\Bbb P}(\sigma(k),z_k)\in\tilde W_k$.}

\medskip

{\bf (i)} Now define $p_5\in\Bbb P$ by setting

\qquad$\tilde I_i'=\tilde I_i\cup\formset{\tilde z(j)(i):j<r}\enskip\forall
\enskip i\in\Bbb N$,

\qquad$p_5(0)=\formset{t:t\in p_4(0),\,C_i(t)\cap\tilde I'_i=\emptyset
\text{ whenever }k_1
\le i<\text{rank}(t)}$,

\qquad$p_5(\xi)=p_4(\xi)$ if $0<\xi<\kappa$.

\noindent Of course we must check that $p_5(0)$, as so defined, belongs to
$\Bbb Q_0\cong\Bbb Q$;  but because $\delta_k(p_4(0))\ge 1$ for $k\ge k_1$,
we have

\centerline{$\text{dp}(\formset{C:t\times C\in p_4(0)})\ge\exp(1/\gamma_k)
\ge 2(\prod_{i<k}2^{n_i})^{k+1}+2k\ge 2\#(\tilde I'_k)$}

\noindent for every $k\ge k_1$, $t\in p_4(0)\cap T_k$, using 1E(a)(vii).
Of course $p_5(0)\cap T_{k_1-1}=p_4(0)\cap T_{k_1-1}=
\formset{\bold v^*(0)} $,
so every element of $p_5(0)$ has a successor in $p_5(0)$, and also

\centerline{$\delta_i(p_5(0))\ge\delta_i(p_4(0))-\gamma_i\ln 2$}

\noindent for $i\ge k_1$.  (Here at last is the key step which depends
on using dp in our measure of `rapidly branching' given in
1E(d).)   Thus $p_5(0)\in\Bbb Q$ and $p_5\in \Bbb P$.
But also

\centerline{$p_5\Vdash_{\Bbb P}\Phi_k\cap\tilde I'_k=\emptyset\enskip
\forall\enskip k\ge k_1$,}

\noindent so that

\centerline{$p_5\Vdash_{\Bbb P}\sigma(k)\notin\tilde I_k,\,
\sigma(k)\ne \tilde z(j)(k)\enskip\forall\enskip j<r,\,(\sigma(k),z_k)
\in\tilde W_k$}

\noindent for $k\ge k_1$;  finally

\centerline{$p_5\Vdash_{\Bbb P}(\sigma,\tilde z)\in\ulcorner R_r\urcorner$,
$\sigma\in\ulcorner R_r\urcorner^{-1}[D']$;} 

\noindent as $p_5\le p$ and $p$, $\sigma$ are arbitrary,

\centerline{$\Bbbone_{\Bbb P}\Vdash_{\Bbb P}\Psi\cap\prod_{k\in\Bbb N}L_k
\subseteq\ulcorner R_r\urcorner^{-1}[D']$,}

\noindent as claimed.

\medskip

{\bf 1S Theorem} For each $r\in\Bbb N$,

$$\eqalign{\Bbbone_{\Bbb P}\Vdash_{\Bbb P}
&\text{ if }D_i\subseteq\ulcorner X\urcorner
\text{ and }D_i\cap\prod_{k\in\Bbb N}
L_k
\ne\emptyset\enskip\forall\enskip\langle L_k\rangle_{k\in\Bbb N}\in\ulcorner
\Cal L\urcorner,\,i<r,\cr
&\text{ then}\enskip\forall\enskip\langle L_k\rangle_{k\in\Bbb N}\in\ulcorner
\Cal L\urcorner\enskip\exists\enskip\langle x_i\rangle_{i<r}\in
\prod_{i<r}(D_i\cap\prod_{k\in\Bbb N}L_k)\cr
&\enskip\enskip\enskip\enskip\text{ such that }
(x_j,x_i)\in \ulcorner R\urcorner\enskip\forall\enskip i<j<r.\cr}$$

\medskip

\noindent {\bf proof}  Induce on $r$.   If $r=0$ the result is trivial.
For the inductive step to $r+1$,  take  
$\Bbb P$-names $\Delta_i$ for  subsets of
$\ulcorner X\urcorner$ such that

\centerline{$\Bbbone_{\Bbb P}\Vdash_{\Bbb P}\Delta_i\cap\prod_{k\in\Bbb N}  
L_k\ne\emptyset\enskip\forall\enskip\langle L_k\rangle_{k\in\Bbb N}
\in\ulcorner\Cal L\urcorner,\,i\le r$.}

\noindent Take a $\Bbb P$-name $\frak L$ for a member of $\ulcorner\Cal L
\urcorner$.   Because members of $\Cal L$ can be coded by simple sequences,
we may suppose that 
$\frak L$ is a $\Bbb P_{\alpha}$-name for some $\alpha<\kappa$ (1Je).   
The inductive hypothesis tells us that

$$\eqalign{\Bbbone_{\Bbb P}\Vdash_{\Bbb P}\,\forall\enskip\langle L_k\rangle
_{k\in\Bbb N}\in\ulcorner\Cal L\urcorner\enskip&\exists\enskip
\langle x_i\rangle_{i<r}\in\prod_{i<r}(\Delta_i\cap\prod_{k\in\Bbb N}L_k)\cr
&\text{such that }(x_j,x_i)\in\ulcorner R\urcorner\text{ if }i<j<r.\cr}$$

\noindent Next, using 1Jc-e,
we can find a $\beta\in\kappa\setminus
\alpha$ and a $\Bbb P_{\beta}$-name $\Delta$ for a subset of $\ulcorner
X^r\urcorner$ such that

\centerline{$\Bbbone_{\Bbb P}\Vdash_{\Bbb P}\Delta\subseteq\prod_{i<r}
\Delta_i$,}

$$\eqalign{\Bbbone_{\Bbb P_{\beta}}\Vdash_{\Bbb P_{\beta}}
(x_j,x_i)\in\ulcorner R\urcorner^{(\beta)}&\text{ whenever }
\langle x_l\rangle_{l<r}\in\Delta,\,i<j<r,\cr
&\Delta\cap(\prod_{k\in\Bbb N}L_k)^r\ne\emptyset\enskip\forall\enskip
\langle L_k\rangle_{k\in\Bbb N}\in\ulcorner\Cal L\urcorner^{(\beta)},\cr}$$

\noindent where we write $\ulcorner\ldots\urcorner^{(\beta)}$ to indicate that
we are interpreting some formula in $V^{\Bbb P_{\beta}}$.

Now we remark that by 1Jb $\Bbb P$ can be regarded, for forcing purposes,
as an iteration $\Bbb P_{\beta}
*\Bbb P'$, where $\Bbb P'$ is a $\Bbb P_{\beta}$-name for a p.o.set
with the same definition, interpreted in $V^{\Bbb P_{\beta}}$, as $\Bbb P$
has in the ground model.
So we may use Lemma 1R in $V^{\Bbb P_{\beta}}$ to say that

\centerline{$\Bbbone_{\Bbb P_{\beta}}\Vdash_{\Bbb P_{\beta}}
\bigl(\Bbbone_{\Bbb P'}\Vdash_{\Bbb P'}\Psi^{(\beta)}\cap\prod\frak L
\subseteq   
\ulcorner R_r\urcorner^{-1}[\Delta\cap(\prod\frak L)^r]\bigr)$,}

\noindent using the notation $\Psi^{(\beta)}$ to indicate which version of
the $\Bbb P$-name $\Psi$ we are trying to use.   Moving to $V^{\Bbb P}$
for a moment, we have $\ulcorner \mu\urcorner\Psi^{(\beta)}=1$ and
$\ulcorner\mu\urcorner(\prod\frak L)>0$, so

\centerline{$\Bbbone_{\Bbb P}\Vdash_{\Bbb P}\exists\enskip l\in\Bbb N,\,
\ulcorner\mu\urcorner(\Psi^{(\beta)}_l\cap\prod\frak L)>0$.}

\noindent Also, of course, every $\Psi_l^{(\beta)}\cap\prod\frak L$
can be regarded (in $V^{\Bbb P}$) as the product of a sequence belonging to  
$\ulcorner\Cal L\urcorner$.
By the original hypothesis on $\Delta_r$,

\centerline{$\Bbbone_{\Bbb P}\Vdash_{\Bbb P} \Delta_r\cap\Psi^{(\beta)}
\cap\prod\frak L\ne\emptyset$.}

\noindent We can therefore find a $\Bbb P$-name $\sigma_r$ for a member of
$\Delta_r\cap\Psi^{(\beta)}\cap\prod\frak L$, and now further 
$\Bbb P_{\beta}$-names
$\sigma_i$, for $i<r$, such that

\centerline{$\Bbbone_{\Bbb P}\Vdash_{\Bbb P}\langle \sigma_i\rangle_{i<r}
\in\Delta\cap(\prod\frak L)^r,
\,(\sigma_r,\langle \sigma_i\rangle_{i<r})\in\ulcorner R_r\urcorner$.}

\noindent But of course we now have

\centerline{$\Bbbone_{\Bbb P}\Vdash_{\Bbb P}
\langle\sigma_i\rangle_{i\le r}\in\prod_{i\le r}(\Delta_i\cap\prod\frak L),\,
(\sigma_j,\sigma_i)\in \ulcorner R\urcorner\enskip\forall\enskip i<j\le r$.}

\noindent 
As $\langle \Delta_i\rangle_{i\le r}$, $\frak L$ are arbitrary, this
shows that the induction proceeds.

\discrpage

\versiondate{6.8.91}

\medskip

\noindent
{\bf 2. Pointwise compact sets of measurable functions}  We turn now to the
questions in analysis which the construction in $\S$1 is designed to solve.
We begin with some definitions and results taken from \TalHD.

\medskip

{\bf 2A Definitions (a)}
Let $(X,\Sigma,\mu)$ be a probability space.   
Write $\Cal L^0=\Cal L^0(\Sigma)
\subseteq\Bbb R^X$ for the set of $\Sigma$-measurable real-valued functions
on $X$.   Let $\frak T_p$ be the topology of pointwise convergence, the usual
product topology, on $\Bbb R^X$.   Let $\frak T_m$ be the (non-Hausdorff,
non-locally-convex) topology of convergence in measure on $\Cal L^0$,
defined by the pseudometric

\centerline{$\rho(f,g)=\int\min(|f(x)-g(x)|,1)\mu(dx)$}

\noindent for $f$, $g\in\Cal L^0$.

\medskip

{\bf (b)} 
A set $A\subseteq\Bbb R^X$ is {\bf stable} if whenever $\alpha<\beta$ in
$\Bbb R$, $E\in\Sigma$ and $\mu E>0$ there are $k$, $l\ge 1$ such that

\quad$\mu^*_{k+l}\formset{(\bold x,\bold y):\bold x\in E^k,\,
\bold y\in E^l,\,\exists\enskip f\in A,\,f(\bold x(i))\le\alpha
\,\&\,f(\bold y(j))\ge\beta\enskip\forall\enskip i<k,\,j<l}$

\centerline{$<(\mu E)^{k+l}$,}

\noindent writing 
$\mu^*_{k+l}$ for the usual product outer measure
on $X^k\times X^l$.   (See \TalHD, 9-1-1.)

\medskip

{\bf 2B Stable sets} Suppose that $(X,\Sigma,\mu)$ is a probability space
and that $A\subseteq\Bbb R^X$ is a stable set.

\medskip

{\bf (a)} If $(X,\Sigma,\mu)$ is complete, then $A\subseteq\Cal L^0(\Sigma)$.
(\TalHD, $\S$9.1.)

\medskip

{\bf (b)}  The $\frak T_p$-closure of $A$ in $\Bbb R^X$ is stable.

\medskip

{\bf (c)} If $A$ is bounded above and below by members of $\Cal L^0$, its
convex hull is stable (\TalHD, 11-2-1).

\medskip

{\bf (d)}  If $A\subseteq\Cal L^0$ (as in (a)), then $\frak T_m\restr A$, the
subspace topology on $A$ induced by $\frak T_m$, is coarser than $\frak T_p
\restr A$.
(\TalHD, 9-5-2.)

\medskip

For more about stable sets, see \TalHD and \TaHG.

\medskip

{\bf 2C Pettis integration} Let $(X,\Sigma,\mu)$ be a probability space and
$B$ a (real) Banach space.

\medskip

{\bf (a)} A function $\phi:X\to B$ is {\bf scalarly measurable}
if $g\phi:X\to\Bbb R$ is $\Sigma$-measurable for every $g\in B^*$, the 
continuous dual of $B$.

\medskip

{\bf (b)} In this case, $\phi$ is {\bf Pettis integrable} if there is a
function $\theta:\Sigma\to B$ such that

\centerline{$\int_E g\phi\, d\mu\text{ exists }=g(\theta E)\enskip\forall\enskip
E\in\Sigma,\,g\in B^*$.}

\medskip

{\bf (c)} If $\phi:X\to B$ is bounded and scalarly measurable, then

\centerline{$A=\formset{g\phi:g\in B^*,\,\|g\|\le 1}\subseteq\Cal L^0$}

\noindent is $\frak T_p$-compact.   In this case $\phi$ is Pettis integrable
iff

\centerline{$f\mapsto\int_Ef:A\to\Bbb R$}

\noindent is $\frak T_p\restr A$-continuous for every $E\in\Sigma$
(\TalHD, 4-2-3).   In
particular (by 2B(d)) $\phi$ is Pettis integrable if $A$ is stable.

\medskip

{\bf 2D The rivals} Write $\mu_L$ for Lebesgue measure on $[0,1]$, and
$\Sigma_L$ for its domain. 
Consider the following two propositions:

{\narrower\smallskip 
\noindent (*) $[0,1]$ is not the union of fewer than $\frak c$ closed negligible
sets;

\medskip

\noindent ($\dagger$) there are sequences $\langle n_k\rangle_{k\in\Bbb N}$,
$\langle W_k\rangle_{k\in\Bbb N}$ such that

$n_k\ge 2^k,\,W_k\subseteq n_k\times n_k,\,\#(W_k)\le 2^{-k}n_k^2
\enskip\forall\enskip k\in\Bbb N$;

taking $X=\prod_{k\in\Bbb N}n_k$, $\mu$ the usual Radon probability on $X$,

\centerline{$R=\formset{(x,y):x,\,y\in X,\,(x(k),y(k))\in W_k\enskip
\forall\enskip k\in\Bbb N,\,\formset{k:x(k)=y(k)}\text{ is finite}}$,}

\noindent then whenever $D\subseteq X$, $\mu^*D=1$ and $r\in\Bbb N$ there
are $x_0,\ldots,x_r\in D$ such that $(x_j,x_i)\in R$ whenever $i<j\le r$.

\smallskip}

\noindent Evidently (*) is a consequence of CH, while
in the language of $\S$1, 
$\Bbbone_{\Bbb P}\Vdash_{\Bbb P}
(\dagger)$, this being a slightly weaker version of Theorem 1S.

Thus both (*) and ($\dagger$) are relatively consistent with ZFC.   
Consequences
of (*) are explored in \TalHD, where it is called Axiom
L;  we list a few of them in 2E below.
Our purpose in this paper is to show that ($\dagger$) leads to a somewhat
different world.

\medskip

{\bf 2E Theorem} Assume (*).   Write $\Cal L^0$ for $\Cal L^0(\Sigma_L)$.

{\bf (a)} If $A\subseteq\Cal L^0$ is separable and compact for
$\frak T_p$, it is stable.

{\bf (b)} If $A\subseteq\Cal L^0$ is separable and compact for $\frak T_p$,
its closed convex hull in $\Bbb R^{[0,1]}$ lies within $\Cal L^0$.

{\bf (c)} If $A\subseteq\Cal L^0$ is separable and compact for $\frak T_p$,
then $\frak T_m\restr A$ is coarser than $\frak T_p\restr A$.

{\bf (d)} If $(Y,\frak S,\text{T},\nu)$ is a separable compact Radon measure
space and $f:[0,1]\times Y\to \Bbb R$ is measurable in the first variable
and continuous in the second, then it is measurable for the (completed)
product measure $\mu_L\times\nu$.

{\bf (e)} If $\langle E_n\rangle_{n\in\Bbb N}$ is a stochastically independent
sequence of measurable subsets of $[0,1]$, with $\lim_{n\to\infty}\mu E_n=0$
but $\sum_{n\in\Bbb N}(\mu E_n)^k=\infty$ for every $k\in\Bbb N$, then
there is an ultrafilter $\Cal F$ on $\Bbb N$ such that

\centerline{$\lim_{n\to\Cal F}E_n=\formset{x:\formset{n:x\in E_n}\in\Cal F}$}

\noindent is non-measurable.

\medskip

\noindent{\bf proof (a)}  See \TalHD, 9-3-1(b).   {\bf (b)} Use (a) and
2B(c).   {\bf (c)} Use (a) and 2B(d).   {\bf (d)} 
Use (a) and \TalHD, 10-2-1.   {\bf (e)} Observe that,
writing $\chi E_n$ for the characteristic function of $E_n$, the set
$\formset{\chi E_n:n\in\Bbb N}$ is not stable, and use (a).

\medskip

{\bf 2F Theorem} Assume ($\dagger$).

{\bf (a)} There is a bounded
Pettis integrable function $\phi:[0,1]\to \ell^{\infty}$ such that
$\formset{g\phi:g\in (\ell^{\infty})^*,\,\|g\|\le 1}$ is not stable in $\Cal L^0
(\Sigma_L)$.

{\bf (b)} There is a separable convex $\frak T_p$-compact subset of
$\Cal L^0(\Sigma_L)$ which is not stable.

\medskip

\noindent{\bf proof} (We write $\ell^{\infty}$ for the Banach space of
bounded real sequences.)  
Take $\langle n_k\rangle_{k\in\Bbb N}$, $\langle W_k\rangle_{k\in\Bbb N}$,
$X$, $\mu$, $R$ from the statement of ($\dagger$).  
Because $([0,1],\mu_L)$ is isomorphic, as measure space,
to $(X,\mu)$,
we may work with $X$
rather than with $[0,1]$.
Write $\Sigma$ for the domain of $\mu$, $\Cal L^0=\Cal L^0(\Sigma)$.

\medskip

{\bf (a)} 
For $k\in\Bbb N$ write

\centerline{$\Cal I_k=\formset{I:I\subseteq n_k,\,
\#(I)\le k,\,(i,j)\notin W_k
\text{ for all distinct }i,\,j\in I}$,}

\noindent For $k\in\Bbb N$, $I\subseteq n_k$ set

\centerline{$H_{kI}=\formset{x:x\in X,\,x(k)\in I}$.}

\noindent Let $A$ be

\centerline{$\formset{\chi H_{kI}:k\in\Bbb N,\,I\in\Cal I_k}$,}

\noindent writing $\chi H:X\to\formset{0,1}$ for the characteristic function
of $H\subseteq X$;  let $Z$ be the $\frak T_p$-closure of $A$ in $\Bbb R^X$.
Because $A$ is uniformly bounded, $Z$ is $\frak T_p$-compact.   
For $E\in\Sigma$ define $f_E:Z\to\Bbb R$ by setting $f_E(h)=\int_E h(x)
\mu(dx)$ for $h\in A$, $f_E(u)=0$ for $u\in Z\setminus A$.   Enumerate
$A$ as $\langle h_m\rangle_{m\in\Bbb N}$, and define $\phi:X\to\ell^{\infty}$,
$\theta:\Sigma\to\ell^{\infty}$ by setting

\centerline{$\phi(x)(m) = h_m(x)\enskip\forall\enskip m\in\Bbb N,\,x\in X$,}

\centerline{$\theta(E)(m) = \int_E h_m(x)\mu(dx)\enskip\forall \enskip
m\in\Bbb N,\,E\in\Sigma$.}

We aim to show 

(i) that $A$ is not stable; 

(ii) that if $\nu$ is a
Radon probability on $A'=Z\setminus A$ then $\int u(x)\,\nu(du)=0$
for $\mu$-almost every $x$; 

(iii) $f_E:Z\to\Bbb R$ is continuous for every $E\in\Sigma$;

(iv) $\theta$ is the indefinite Pettis integral of $\phi$, so that
$\phi$ is Pettis integrable;

(v) $K=\formset{g\phi : g\in C(Z)^*,\,\|g\|\le 1}$ includes $A$ so
is not stable.

\medskip

{\bf ad (i)} Suppose that $k$, $l\ge 1$.  Take any $m\ge l$. 
Set

$$\eqalign{G&=\formset{\bold y:\bold y\in X^l,\,\exists\enskip I\in\Cal I_m,
\,\bold y(j)\in H_{mI}\enskip\forall\enskip j<l}\cr
&\supseteq
\formset{\bold y:\bold y\in X^l,\,(\bold y(i)(m),\bold y(j)(m))\notin W_m
\text{ for distinct }i,\,j<l}.\cr}$$

\noindent Because $\#(W_m)\le 2^{-m}n_m^2$, $\mu_l G\ge (1-2^{-m})^{l^2}$.
If $\bold y\in G$, set $I=\formset{\bold y(j)(m):j<l}\in\Cal I_m$;
then

\centerline{$\mu_k\formset{\bold x:\bold x\in X^k,\,\bold x(i)(m)\notin I\enskip
\forall\enskip i<k}\ge(1-n_m^{-1}l)^k$.}

\noindent So we conclude that

\centerline{$\mu_{k+l}\formset{(\bold x,\bold y):\bold x\in X^k,\,
\bold y\in X^l,\,\exists\enskip f\in A,\,f(\bold x(i))=0\enskip\forall\enskip
i<k,\,f(\bold y(j))=1\enskip\forall\enskip j<l}$}

\centerline{$\ge (1-2^{-m})^{l^2}(1-n_m^{-1}l)^k$}

\noindent (by Fubini's theorem).   Because $k$, $l$ and $m$ are arbitrary,
$A$ cannot be stable.

\medskip

{\bf ad (ii)} 
Because each $\Cal I_m$ is finite, any member of $A'$
must be of the form $\chi E$ where $E\subseteq X$ and

\centerline{$x\in E,\,x'\in X,\,\formset{k:x(k)\ne x'(k)}\text{ is finite }
\Rightarrow x'\in E$.}

\noindent Note also that if $x$, $y\in E$ then $(x,y)\notin R$;  because
either $x(k)=y(k)$ for infinitely many $k$, or there are $k$, $I$ such that
$x(k)\ne y(k)$, $I\in\Cal I_k$ and $x$, $y$ both belong to $H_{kI}$, in
which case $(x(k),y(k))\notin W_k$.

Now let $\nu$ be a Radon probability
on $A'$, and set $w(x)=\int u(x)\,\nu(du)$ for each $x\in X$, 
so that $w$ belongs to the closed convex hull of $A'$ in
$\Bbb R^X$.  If $x$, $x'$ are two members of $X$ differing on only finitely
many coordinates, then $u(x)=u(x')$ for every $u\in A'$;  consequently
$w(x)=w(x')$.   Also $0\le w(x)\le 1$ for every $x\in X$.

Take $\delta>0$ and set $D=\formset{x:w(x)\ge\delta}$.   By the zero-one
law, $\mu^*D$ must be either $0$ or $1$.   Suppose, if possible, that
$\mu^*D=1$.   Let $r\in\Bbb N$ be such that $r\delta\ge 1$.  
By ($\dagger$), there are $x_0,\ldots,x_r\in D$
such that $(x_j,x_i)\in R$ for $i<j\le r$.   But
in this case $\sum_{i\le r}u(x_i)\le 1$ for every $u\in A'$, while
$\sum_{i\le r}w(x_i)\ge(r+1)\delta>1$, and $w$ cannot belong to the
closed convex hull of $A'$.

Accordingly $\mu^*D$ must be $0$.   As $\delta$ is arbitrary, $w
=0$ a.e.

\medskip

{\bf ad (iii)} Because $f_E(\chi H_{kI})\le kn_k^{-1}$ for every
$I\in\Cal I_k$, $\lim_{m\to\infty}f_E(h_m)=0$ and $f_E$ is continuous.

\medskip

{\bf ad (iv)} We need to show that

\centerline{$\int_Eg(\phi(x))\,\mu(dx)\text{ exists }=
g(\theta(E))\enskip\forall\enskip
g\in (\ell^{\infty})^*,\,E\in\Sigma$.}

\noindent It is enough to consider positive linear functionals $g$ of
norm 1.   For any such $g$ we have a Radon probability $\nu$ on $Z$
such that 

\centerline{$g(\langle f(h_m)\rangle_{m\in\Bbb N})
=\int_Z f(u)\,\nu(du)$ for every $f\in C(Z)$,}

\noindent using the Riesz representation of positive linear functionals
on $C(Z)$.
Set $\epsilon_m=\nu\formset{h_m}$, $\epsilon=1-\sum_{m\in\Bbb N}
\epsilon_m=\nu A'$.   Then we can find a Radon probability $\nu'$
on $A'$ such that

\centerline{$g(\langle f(h_m)\rangle_{m\in\Bbb N})
=\sum_{m\in\Bbb N}\epsilon_mf(h_m)+\epsilon\int_{A'}f(u)
\,\nu'(du)$}

\noindent for every $f\in C(Z)$.   Now an easy calculation (using
(ii)) shows that

\centerline{$g(\theta(E))=g(\langle f_E(h_m)\rangle_{m\in\Bbb N})
=\sum_{m\in\Bbb N}\epsilon_m\int_E
h_m(x)\,\mu(dx)=\int_Eg(\phi(x))\,\mu(dx)$}

\noindent for every $E\in\Sigma$.

\medskip

{\bf ad (v)} If $m\in\Bbb N$ then $h_m=e_m\phi\in K$ where 
$e_m\in(\ell^{\infty})^*$
is defined by setting $e_m(z)=z(m)$ for every $z\in\ell^{\infty}$.
This completes the proof.

\medskip

{\bf (b)} 
The unit ball of $(\ell^{\infty})^*$ is w$^*$-separable and its continuous image
$K\subseteq\Cal L^0$ is separable;  so $K$ witnesses the truth of (b).

\medskip

{\bf 2G Further properties of the model} Returning to 1R/1S, we see that the
model of $\S$1 has some further striking characteristics closely allied to, but
not obviously derivable from, ($\dagger$).   Consider for instance

{\narrower\smallskip
\noindent ($\ddagger$) there is a closed negligible set $Q\subseteq [0,1]^2$
such that whenever $D\subseteq [0,1]$ and $\mu_L^*D=1$ then $\mu_LQ^{-1}[D]=1$;

\medskip

\noindent ($\ddagger$)$'$ there is a negligible set $Q'\subseteq [0,1]^2$
such that whenever $C$, $D\subseteq [0,1]$ and $(C\times D)\cap Q'=\emptyset$
then one of $C$, $D$ is negligible.

\smallskip}

\noindent Then $\Bbbone_{\Bbb P}\Vdash_{\Bbb P}(\ddagger)$.  For start by taking
$Q_1$ to be

\centerline{$\formset{(x,y):x,\,y\in X,\,(x(k),y(k))\in W_k\enskip
\forall\enskip k\in\Bbb N}$,}

\noindent the closure of $R$ in $X\times X$.   Then the argument for 1S shows
that

$$\eqalign{\Bbbone_{\Bbb P}\Vdash_{\Bbb P}\text{ if }&D\subseteq\ulcorner
X\urcorner\text{ and }D\cap\prod_{k\in\Bbb N}L_k\ne\emptyset\enskip
\forall\enskip\langle L_k\rangle_{k\in\Bbb N}\in\ulcorner\frak L\urcorner\cr
&\text{then }\exists\enskip
\beta<\kappa\text{ such that }\Psi^{(\beta)}\subseteq
\ulcorner Q_1\urcorner^{-1}[D].\cr}$$

\noindent Consequently

\centerline{$\Bbbone_{\Bbb P}\Vdash_{\Bbb P}\text{if }D\subseteq\ulcorner
X\urcorner\text{ and }\ulcorner\mu\urcorner^*D=1\text{ then }
\ulcorner\mu\urcorner(\ulcorner Q_1\urcorner^{-1}[D])=1$.}

\noindent Accordingly we have in $V^{\Bbb P}$ the version of $(\ddagger)$
in which $([0,1],\mu_L)$ is replaced by $(X,\mu)$. 
However there is now
a continuous inverse-measure-preserving function $f:X\to[0,1]$, and taking

\centerline{$Q=\formset{(f(x),f(y)):(x,y)\in Q_1}$}

\noindent we obtain $(\ddagger)$ itself.   Evidently $(\ddagger)$
implies $(\ddagger)'$, taking $Q'$ to be

\centerline{$\formset{(x+q,y+q'):(x,y)\in Q,\,q,\,q'\text{ are rational}}
\cap[0,1]^2$.}

Of course $(^*)$ and $(\ddagger)$
are mutually incompatible (the argument for 2E(a) from $(^*)$, greatly simplified, demolishes $(\ddagger)$ also).    The weaker form $(\ddagger)'$ is
incompatible with CH or MA,  but not with $(^*)$, both $(\ddagger)'$
and $(^*)$ being true in Cohen's original model of not-CH (see \FrpHI).

\medskip

{\bf 2H Problems} The remarkable results quoted in 
2E depend on the identification
of separable relatively pointwise compact sets with stable sets (`Axiom F'
of \TalHD).   In models satisfying ($\dagger$), this identification breaks
down.   But our analysis does not seem to touch any of 2E(b)-(e).   We
therefore spell out the obvious problems still outstanding.   Write
$\Cal L^0$ for $\Cal L^0(\Sigma_L)$. 

\medskip

{\bf (a)} Is it relatively consistent with ZFC to suppose that there is
a separable $\frak T_p$-compact set $A\subseteq\Cal L^0$ such that the
closed convex hull of $A$ in $\Bbb R^{[0,1]}$ does not lie within $\Cal L^0$?

\medskip

{\bf (b)} Is it relatively consistent with ZFC to suppose that there is a
separable $\frak T_p$-compact set $A\subseteq\Cal L^0$ such that 
$\frak T_m\restr A$ is not coarser than $\frak T_p\restr A$?   Does it
make a difference if $A$ is assumed to be convex?   (This question seems
first to have been raised by J.Bourgain and F.Delbaen.)

\medskip

{\bf (c)} Is it relatively consistent with ZFC to suppose that there are
a separable compact Radon measure space $(Y,\frak S,\text{T},\nu)$ and
a function $f:[0,1]\times Y\to\Bbb R$ which is measurable in the first
variable, continuous in the second variable, but not jointly measurable
for $\mu_L\times \nu$?

\medskip

{\bf (d)} Is it relatively consistent with ZFC to suppose that there is a
stochastically independent sequence $\langle E_n\rangle_{n\in\Bbb N}$
in $\Sigma_L$ such that $\sum_{n\in\Bbb N}(\mu_L E_n)^k=\infty$ for every
$k\in\Bbb N$, but $\mu_L(\lim_{n\to\Cal F}E_n)=0$ for every non-principal
ultrafilter $\Cal F$ on $\Bbb N$?   (This question is essentially due to
W.Moran;  see also \TalHD, 9-1-4 for another version.)

\medskip

Here we note only that a positive answer to (a) would imply the same answer
to (c), and that the word `separable' in (a)-(c) is necessary, as is shown
by examples 3-2-3 and 10-1-1 in \TalHD.

\bigskip

{\bf \noindent References}

\BJSpHI\ T.Bartoszy\'nski, H.Judah \& S.Shelah, `The Cicho\'n diagram',
preprint, 1989 (MSRI 00626-90).

\BaHD\ J.E.Baumgartner, `Applications of the proper forcing axiom',
pp.~913-959 in \KVHD.



\FrpHI\ H.Friedman, `Rectangle inclusion problems', Note of 9 October 1989.

\KuHJ\ K.Kunen, {\it Set Theory.}   North-Holland, 1980.

\KVHD\ K.Kunen \& J.E.Vaughan (eds.), {\it Handbook of Set-Theoretic Topology.}
North-Holland, 1984.

\ShHB\ S.Shelah, {\it Proper Forcing.}   
Springer, 1982 (Lecture Notes in Mathematics 940).

\ShCBF\ S.Shelah, `Vive la diff\'erence!', submitted for the proceedings of
the set theory conference at MSRI, October 1989;  notes of July 1987, preprint
of October 1989;  abbreviation `ShCBF'.

\SpHG\ J.Spencer, {\it Ten Lectures on the Probabilistic Method,}  
S.I.A.M., 1987.

\TalHD\  M.Talagrand, {\it Pettis integral and measure theory}. Mem. Amer. Math.
Soc. 307 (1984).

\TaHG\ M.Talagrand, `The Glivenko-Cantelli problem', Ann. of Probability
15 (1987) 837-870.

\bigskip

\noindent{\bf Acknowledgements} Part of the work of this paper was done
while the authors were visiting the M.S.R.I., Berkeley;  we should like to
thank the Institute for its support.  The first author was partially
supported by the Fund for Basic Research of the Israel Academy of 
Sciences.   The second author was partially
supported by grants GR/F/70730 and GR/F/31656 from the U.K. Science and
Engineering Research Council.
We are most grateful to M.Burke for carefully checking the manuscript.

\bigskip

\discrversionC{}{To appear in J.S.L.}

\end